\documentclass[10pt]{article}%
\usepackage{latexsym,amsmath,amssymb,graphics,amscd}
\usepackage{amsmath}
\usepackage{amsfonts}
\usepackage{amssymb}
\usepackage{graphicx}%
\newcommand{\bfi}{\bfseries\itshape}
\addtocounter{footnote}{1}
\makeatletter
\@addtoreset{figure}{section}
\def\thefigure{\thesection.\@arabic\c@figure}
\def\fps@figure{h,t}
\@addtoreset{table}{bsection}
\def\thetable{\thesection.\@arabic\c@table}
\def\fps@table{h, t}
\@addtoreset{equation}{section}

\makeatother
\newtheorem{theorem}{Theorem}[section]
\newtheorem{definition}[theorem]{Definition}

\newtheorem{proposition}[theorem]{Proposition}

\newtheorem{example}[theorem]{Example}

\newsavebox{\savepar}

\makeatletter

\begin{document}

\title{\textbf{The stochastic Hamilton-Jacobi equation}}
\author{Joan-Andreu L\'{a}zaro-Cam\'{\i}$^{1}$ and Juan-Pablo Ortega$^{2}$}
\date{}
\maketitle

\begin{abstract}
We extend some aspects of the Hamilton-Jacobi theory to the category of  stochastic Hamiltonian dynamical systems. More specifically, we show that the stochastic action satisfies the Hamilton-Jacobi equation when, as in the classical situation, it is written as a function of the configuration space using a regular Lagrangian submanifold. Additionally, we will use a variation of the Hamilton-Jacobi equation to characterize the generating functions of one-parameter groups of symplectomorphisms that allow to rewrite a given stochastic Hamiltonian system in a form whose solutions are very easy to find; this result recovers in the stochastic context the classical solution method by {\it reduction to the equilibrium} of a Hamiltonian system.

\end{abstract}

\medskip

\footnotesize
\noindent {\bf Keywords}: stochastic differential equation, Hamiltonian stochastic differential equation, Hamilton-Jacobi equation.
\normalsize

\makeatletter
\addtocounter{footnote}{1} \footnotetext{Departamento de F\'{\i}sica
Te\'{o}rica. Universidad de Zaragoza. Pedro Cerbuna, 12. E-50009 Zaragoza.
Spain. {\texttt{lazaro@unizar.es}}} \addtocounter{footnote}{1}
\footnotetext{Centre National de la Recherche Scientifique, D\'{e}partement de
Math\'{e}matiques de Besan\c{c}on, Universit\'{e} de Franche-Comt\'{e}, UFR
des Sciences et Techniques. 16, route de Gray. F-25030 Besan\c{c}on cedex.
France. {\texttt{Juan-Pablo.Ortega@univ-fcomte.fr} }} \makeatother

\section{Introduction}

Hamiltonian diffusions  were introduced and  studied by Bismut in the monograph~\cite{Bismut}. These systems were generalized in~\cite{paper hamiltoniano} to accommodate arbitrary Poisson manifolds as phase spaces and general continuous semimartingales as forcing noises. In that paper it was also shown that, when the phase space is an exact symplectic manifold, the stochastic Hamilton equations are fully characterized by a variational principle that generalizes the classical Hamilton's Principle.  This circle of problems has also been treated  in~\cite{bou rabee 1} in the development of stochastic variational numerical integrators.

Hamilton-Jacobi theory is an important part of classical mechanics that provides a characterization of the generating functions  of certain time-dependent canonical transformations that put a given Hamiltonian system in such a form that his solutions are extremely easy to find; this is the so called solution by {\it reduction to the equilibrium}. In this respect, the fact that the classical action satisfies the Hamilton-Jacobi equation is a very relevant result. Hamilton-Jacobi theory also plays a fundamental role in the study of the quantum-classical relationship, in  integrable systems, or in the development of structure preserving numerical integrators. For all these reasons it is desirable to have at hand similar tools in the stochastic Hamiltonian context; this is the main goal of this work. The Hamilton-Jacobi equation was already studied by Bismut~\cite{Bismut} in the context of Hamiltonian diffusions and, as we will see, most of the ideas in that piece of work are still valid at our degree of generality; at some level, this paper can be seen as a completion of Bismut's work in which complete proofs are provided and where the results have been adapted to our framework using a more modern geometric language; this makes them more palatable to a growing community interested both in geometric mechanics and in stochastics.

The paper starts with a brief presentation in Section~\ref{Stochastic Hamiltonian dynamical systems and the stochastic action} of some basic facts about stochastic Hamiltonian systems and, more importantly, with the introduction of the stochastic Hamiltonian action. Section~\ref{The stochastic action on Lagrangian submanifolds and the Hamilton-Jacobi equation} is dedicated to showing that the stochastic action satisfies a generalized version  of the Hamilton-Jacobi equation when written as a function of the configuration space using a Lagrangian submanifold (see Theorem~\ref{theorem Hamilton-Jacobi}).  As an application of the results in this section we show in Example~\ref{feynman kac action} how the exponential of the expectation  of the so called projected stochastic action can be used to construct solutions of the heat equation corrected with a potential, in a way that strongly resembles the Feynman-Kac formula. 

The paper concludes with a section on the relation between the solutions of the Hamilton-Jacobi equation and the generating functions of time dependent diffeomorphisms that allow the integration of the Hamiltonian stochastic differential equation in question in an easy manner. The natural framework for carrying this out is that of time-dependent Hamiltonian systems; that is why we have included a subsection that briefly recalls the classical theory of non-autonomous Hamiltonian systems and presents it in a form that is suitable for generalization in the stochastic context. Some of the statements in this section are either inspired or are a direct generalization of analogous results in~\cite{Bismut}; we have nevertheless included them in order to have a complete and self-contained presentation of the theory.

\medskip

\noindent {\bf Acknowledgements:} the authors thank the hospitality of the Centre de Recerca Matem\`atica of the Universitat Aut\`onoma de Barcelona during the program ``Equivariant Problems in Symplectic Geometry", organized by Eva Miranda. This paper was written while the authors took part in that program. J.-A. L.-C. acknowledges support
from the Spanish Ministerio de Educaci\'{o}n y Ciencia grant number
BES-2004-4914. He also acknowledges partial support from MEC grant
BFM2006-10531 and Gobierno de Arag\'{o}n grant DGA-grupos consolidados
225-206.

\section{Stochastic Hamiltonian dynamical systems and the stochastic action}
\label{Stochastic Hamiltonian dynamical systems and the stochastic action}

In this section we recall the basic facts about stochastic
Hamiltonian dynamical systems and we fix the notation that we are going to use
throughout the paper. The content of this section is extracted from
\cite{paper hamiltoniano}, which the reader is encouraged to check with for a
more comprehensive approach to stochastic Hamiltonian systems and their most
relevant properties. Although such systems may be considered on any Poisson
manifold $\left(  M,\left\{  \cdot,\cdot\right\}  \right)  $, we are only
going to deal with the exact symplectic case, that is, stochastic Hamiltonian
systems defined on exact symplectic manifolds. 

The reason for such a restriction is
that in that case there exists a {\bfi  stochastic action}
defined on the set of manifold valued semimartingales such that the stochastic
Hamilton equations can be characterized using a {\bfi  critical action
principle}. In other words, taking suitable variations on the space of
manifold valued semimartingales, a given semimartingale is critical for the
stochastic action if and only if it is a solution of the stochastic Hamiltonian
equations. This generalizes to the stochastic context the way in which the
Hamilton equations are characterized in the classical deterministic contex. Like in
that framework, the stochastic action plays a prominent role in the
description of a Hamiltonian system and, as we will see later on, it also satisfies a stochastic analog of the
{\bfi  Hamilton-Jacobi equation}. 

Let $\left(  M,\omega\right)  $ be a symplectic manifold. Using the
nondegeneracy of the symplectic form $\omega$, one can associate to each function
$h\in C^{\infty}(M)$ a vector field $X_{h}\in\mathfrak{X}\left(  M\right)  $
characterized by the equality
\begin{equation}
{\bf i}  _{X_{h}}\omega=\mathbf{d}h \label{HJ 2}%
\end{equation}
We will say that $X_{h}$ is the {\bfseries\itshape Hamiltonian vector field}
associated to the Hamiltonian function $h$. The expression (\ref{HJ 2}) is
referred to as the {\bfseries\itshape Hamilton equations}. Given $h\in
C^{\infty}\left(  M\right)  $, solving the associated Hamilton equations amounts to finding
the integral curves of $X_{h}$. 

We now introduce the stochastic generalization of~(\ref{HJ 2}). Let $V$ be a real finite dimensional vector space and let $f:M\rightarrow
V$ be a differentiable function taking values in $V$. We define the
differential $\mathbf{d}f:TM\rightarrow V$ as the map given by $\mathbf{d}%
f=p_{2}\circ Tf$, where $Tf:TM\rightarrow TV=V\times V$ is the tangent map of
$f$ and $p_{2}:V\times V\rightarrow V$ is the projection onto the second
factor. If $V=\mathbb{R}$ this definition coincides with the usual
differential. If $\{e_{1},...,e_{r}\}$ is a basis of $V$ and $f=\sum_{i=1}%
^{r}f^{i}e_{i}$, $f^{i}\in C^{\infty}\left(  M\right)  $, then $\mathbf{d}%
f=\sum_{i=1}^{r}\mathbf{d}f^{i}e_{i}$. A
{\bfseries\itshape
Stratonovich operator} from $V$ to $M$ is a family $\{S(v,z)\}_{v\in V,z\in
M}$ such that $S\left(  v,z\right)  :T_{v}V\rightarrow T_{z}M$ is a linear
mapping that depends smoothly on its two entries. The adjoint of $S\left(
v,z\right)  $ is usually denoted by $S^{\ast}\left(  v,z\right)  :T_{z}^{\ast
}M\rightarrow T_{v}^{\ast}V$, $v\in V$ and $z\in M$. Given a
smooth function $h:M\rightarrow V$, $h=\sum_{i=1}^{r}h^{i}e_{i}$, we define
the associated Hamiltonian Stratonovich operator $H\left(  v,z\right)  :T_{v}V\rightarrow T_{z}M$ by
\begin{equation}
H\left(  v,z\right)  (u)=\sum_{i=1}^{r}\left\langle e^{i},u\right\rangle
X_{h_{i}}(z)\text{, \quad}u\in T_{v}V. \label{HJ 3}%
\end{equation}
In this expression, $\{e^{1},...,e^{r}\}$ is the dual basis of $\{e_{1}%
,...,e_{r}\}$. It is easy to check that the adjoint $H^{\ast}(v,z):T_{z}%
^{\ast}M\rightarrow T_{v}^{\ast}V$ of $H\left(  v,z\right)  $ is $H^{\ast
}(v,z)\left(  \alpha_{z}\right)  =-\mathbf{d}h\left(  z\right)  \cdot
\omega^{\#}(z)\left(  \alpha_{z}\right)  $, where $\omega^\sharp (z): T ^\ast _zM \rightarrow T _zM $ is the isomorphism induced by the symplectic form  $\omega  $.

The key to generalizing the Hamilton equations (\ref{HJ 2}) to the stochastic context consists of realizing that they may be restated by saying that a smooth curve $\gamma:[0,T]\rightarrow M$
is an integral curve of the Hamiltonian vector field $X_{h}$, $h\in C^{\infty
}\left(  M\right)  $, if and only if for any $\alpha\in\Omega\left(  M\right)
$ and for any $t\in\lbrack0,T]$%
\begin{equation}
\int_{\left.  \gamma\right\vert _{[0,T]}}\alpha=-\int_{0}^{t}\mathbf{d}%
h\left(  \omega^{\#}(\alpha)\right)  \circ\gamma\left(  s\right)  ds
\label{HJ 28}%
\end{equation}
(\cite[Proposition 2.1]{paper hamiltoniano}). Using this observation
we will define the stochastic
Hamilton equations by specifying the result of integrating an arbitrary one
form $\alpha\in\Omega\left(  M\right)  $ along them. More especifically, let $\left(  \Omega,\{\mathcal{F}_{t}\}_{t\in\mathbb{R}_{+}},P\right)
$ be a filtered probability space,
$X:\mathbb{R}_{+}\times\Omega\rightarrow V$ a continuous semimartingale (that is,
the paths $X_{t}\left(  \cdot\right)  :\Omega\rightarrow V$ are continuous
a.s. for any $t\in\mathbb{R}_{+}$) that takes values on the vector space $V$
with $X_{0}=0$, and $h:M\rightarrow V^{\ast}$ a smooth function. We will say
that a $M$-valued semimartingale $\Gamma:\mathbb{R}_{+}\times\Omega
\rightarrow M$ is a solution of the stochastic Hamilton equations with
{\bfseries\itshape stochastic component} $X$ and {\bfseries\itshape
Hamiltonian function} $h$ if for any $\alpha\in\Omega\left(  M\right)  $%
\begin{equation}
\int\left\langle \alpha,\delta\Gamma\right\rangle =-\int\left\langle
\mathbf{d}h\left(  \omega^{\#}(\alpha)\right)  ,\delta X\right\rangle ,
\label{HJ 4}%
\end{equation}
where the symbol $\delta$ denotes Stratonovich integration.
In other words, the processes that solve the stochastic Hamilton equations are no longer
driven by the deterministic \textit{time} $t$ but by the stochastic
\textit{noise} $X$. It can be shown (\cite[Theorem 7.21]{emery}) that given a
semimartingale $X$ in $V$, a $\mathcal{F}_{0}$-measurable random variable
$\Gamma_{0}$, there are a maximal stopping time $\zeta>0$ and a continuous
solution $\Gamma$ of (\ref{HJ 4}) with initial condition $\Gamma_{0}$ defined
on the set $\{\left(  t,\eta\right)  \in\mathbb{R}_{+}\times\Omega
~|~t\in\lbrack0,\zeta(\eta))\}$. If $\zeta$ is finite, then $\Gamma$ explodes
at time $\zeta$, that is, the path $\Gamma_{t}$ with $t\in\lbrack0,\zeta)$ is
not contained in any compact subset of $M$. For the sake of simplicity, the
stochastic Hamilton equations (\ref{HJ 4}) will be symbolically denoted using the Stratonovich operator $H$ as
\begin{equation}
\delta\Gamma=H\left(  X,\Gamma\right)  \delta X. \label{HJ 5}%
\end{equation}

We now  state some
of the basic properties of the flow defined by~(\ref{HJ 4}). Let $\varphi\left(
\cdot,z\right)  :[0,\zeta(z))\subseteq\mathbb{R}_{+}\times\Omega\rightarrow M$
denote the unique solution of (\ref{HJ 5}) with initial condition $\Gamma
_{0}=z\in M$ a.s.. The map $\varphi$ will be referred to as the
{\bfseries\itshape stochastic flow} associated to (\ref{HJ 5}). For any
$\left(  t,\eta\right)  \in\mathbb{R}_{+}\times\Omega$, let $\mathbb{D}%
_{t}\left(  \eta\right)  =\{z\in M~|~\zeta\left(  z,\eta\right)  >t\}$.
Observe that $\mathbb{D}_{t}\left(  \eta\right)  \subseteq\mathbb{D}%
_{s}\left(  \eta\right)  $ if $s\leq t$. By \cite[Lemma 4.8.3]{Kunita book}
$\mathbb{D}_{t}\left(  \eta\right)  $ is an open set for any $t\in
\mathbb{R}_{+}$ a.s. and%
\[%
\begin{array}
[c]{rrl}%
\varphi_{t}\left(  \eta\right)  :\mathbb{D}_{t}\left(  \eta\right)  &
\longrightarrow & M\\
z & \longmapsto & \varphi_{t}\left(  z,\eta\right)
\end{array}
\]
is a continuously differentiable diffeomorphism (\cite[Theorem 4.8.4]{Kunita
book}). Additionally,
\[%
\begin{array}
[c]{rrl}%
\varphi\left(  \eta\right)  :[0,t]\times\mathbb{D}_{t}\left(  \eta\right)  &
\longrightarrow & M\\
z & \longmapsto & \varphi_{t}\left(  z,\eta\right)
\end{array}
\]
is continuous and its partial derivatives with respect to $z\in\mathbb{D}%
_{t}\left(  \eta\right)  $ are also continuous on $[0,t]\times\mathbb{D}%
_{t}\left(  \eta\right)  $. The local version of these results, that is, the
case $M=\mathbb{R}^{2n}$, can be also found in \cite[Chapter V Theorem
39]{protter book}. Furthermore, the stochastic flow $\varphi$ acts naturally on tensor fields and in particular on differential
forms. Hence, by \cite[Theorem 3.3]{kunita seminaire} and \cite[Section
4.9]{Kunita book}, if $\alpha\in\Omega^{k}\left(  M\right)  $ is a $k$-form,
$k\in\mathbb{N}$, then%
\begin{equation}
\varphi_{t}\left(  \eta\right)  ^{\ast}\alpha=\alpha+\sum_{i=1}^{r}\left(
\int_{0}^{t}\varphi_{s}^{\ast}\left(  \pounds _{X_{h_{i}}}\alpha\right)
\delta X_{s}^{i}\right)  \left(  \eta\right)  \label{HJ 6}%
\end{equation}
on $\mathbb{D}_{t}\left(  \eta\right)  $, $\left(  t,\eta\right)
\in\mathbb{R}_{+}\times\Omega$. In particular, if $\alpha=\omega$ is the
symplectic form, then $\pounds _{X_{h_{i}}}\omega=0$ for any $i=1,...,r$ and
$\varphi^{\ast}\omega=\omega$ which is the stochastic version of the
Liouville's Theorem (see \cite[Theorem 2.1]{paper hamiltoniano}).

We conclude this brief summary on stochastic Hamiltonian systems by
introducing the stochastic action and presenting how it characterizes the
solutions of the Hamilton equations. As we already said,  the stochastic action is only naturally defined  for stochastic
Hamiltonian systems occurring on {\bfi  exact} symplectic manifolds: let
$\left(  M,\omega=-\mathbf{d}\theta\right)  $ be an exact symplectic manifold,
$X:\mathbb{R}_{+}\times\Omega\rightarrow V$ a semimartingale taking values on
the vector space $V$, and $h:M\rightarrow V^{\ast}$a Hamiltonian function. We
denote by $\mathcal{S}\left(  M\right)  $ and $\mathcal{S}\left(
\mathbb{R}\right)  $ the sets of $M$ and real valued semimartingales,
respectively. We define the {\bfseries\itshape stochastic action} associated
to $h$ as the map $S:\mathcal{S}\left(  M\right)  \rightarrow\mathcal{S}%
\left(  \mathbb{R}\right)  $ given by
\[
S\left(  \Gamma\right)  :=\int\left\langle \theta,\delta\Gamma\right\rangle
-\int\langle\hat{h}(\Gamma),\delta X\rangle
\]
where in the previous expression, the stochastic differential form $\hat
{h}\left(  \Gamma\right)  :\mathbb{R}_{+}\times\Omega\rightarrow V\times
V^{\ast}$ over $X$ is given by $\hat{h}\left(  \Gamma\right)  \left(
t,\omega\right)  :=\left(  X_{t}\left(  \omega\right)  ,h\left(  \Gamma
_{t}\left(  \omega\right)  \right)  \right)  $.

Given a $M$-valued semimartingale $\Gamma$ and $s_{0}>0$, we say that the map
$\Sigma:\left(  -s_{0},s_{0}\right)  \times\mathbb{R}_{+}\times\Omega
\rightarrow M$ is a {\bfseries\itshape pathwise variation} of $\Gamma$
whenever $\Sigma_{t}^{s_{0}=0}=\Gamma_{t}$ a.s.. We say that the pathwise
variation $\Sigma$ of $\Gamma$ {\bfseries\itshape converges uniformly} to
$\Gamma$ whenever the following properties are satisfied (\cite[Definition
4.4]{paper hamiltoniano}):

\begin{enumerate}
\item[{\bf  (i)}] For any $f\in C^{\infty}\left(  M\right)  $, $f(\Sigma
^{s})\rightarrow f(\Gamma)$ {\bfseries\itshape uniformly in compacts in
probability} (abbreviated \textit{in ucp}) as $s\rightarrow0$. That is, for
any $\varepsilon>0$ and any $t\in\mathbb{R}_{+}$,
\[
P\left(  \left\{  \sup_{0\leq u\leq t}\left\vert f\left(  \Sigma_{u}%
^{s}\right)  -f\left(  \Gamma_{u}\right)  \right\vert >\varepsilon\right\}
\right)  \underset{s\rightarrow0}{\longrightarrow}0.
\]

\item[{\bf (ii)}] There exists a process $Y:\mathbb{R}_{+}\times\Omega\rightarrow
TM$ over $\Gamma$ such that, for any $f\in C^{\infty}\left(  M\right)  $, the
Stratonovich integral $\int Y[f]\delta X$ exists for any continuous real
semimartingale $X$ (this is for instance guaranteed if $Y$ is a
semimartingale) and, additionally, the increments $\left.  \left(  f\left(
\Sigma^{s}\right)  -f\left(  \Gamma\right)  \right)  \right/  s$ converge in
ucp to $Y[f]$ as $s\rightarrow0$. We will call such a $Y$ the
{\bfseries\itshape infinitesimal generator} of $\Sigma$.
\end{enumerate}

We will say that $\Sigma$ (respectively $Y$) is {\bfseries\itshape bounded}
when its image lies in a compact set of $M$ (respectively $TM$). It can be
shown that, given a $M$-valued semimartingale $\Gamma$, a compact set
$K\subseteq M$, and a bounded process $Y:\mathbb{R}_{+}\times\Omega\rightarrow
TM$ over $\Gamma^{\tau_{K}}$ (the process $\Gamma$ stopped at the first exit
time $\tau_{K}$ of $\Gamma$ from $K$) such that $\int Y[f]\delta X$ exists for
any continuous real semimartingale $X$ and any $f\in C^{\infty}\left(
M\right)  $, there exists a bounded pathwise variation $\Sigma$ that converges
uniformly to $\Gamma^{\tau_{K}}$ whose infinitesimal generator is $Y$
(\cite[Proposition 4.2]{paper hamiltoniano}). Using these elements, a variational characterization of the stochastic Hamilton equations can be given (\cite[Theorem 4.2]{paper hamiltoniano}): the semimartingale
$\Gamma$ satisfies the stochastic Hamilton equations (\ref{HJ 4}) with
initial condition $\Gamma_{t=0}=m_{0}\in M$ a.s. up to time $\tau_{K}$ if and
only if, for any bounded pathwise variation $\Sigma:\left(  s_{0}%
,s_{0}\right)  \times\mathbb{R}_{+}\times\Omega\rightarrow M$ with bounded
infinitesimal generator which converges uniformly to $\Gamma^{\tau_{K}}$ and
such that $\Sigma_{0}^{s}=m_{0}$ and $\Sigma_{\tau_{K}}^{s}=\Gamma_{\tau_{k}}$
a.s., $s\in\left(  -s_{0},s_{0}\right)  $,%
\[
\left[  \left.  \frac{d}{ds}\right\vert _{s=0}S\left(  \Sigma^{s}\right)
\right]  _{\tau_{K}}=0\text{ \ a.s.}.
\]

\section{The stochastic action on Lagrangian submanifolds and the Hamilton-Jacobi equation}
\label{The stochastic action on Lagrangian submanifolds and the Hamilton-Jacobi equation}

It is a classical result in mechanics that the action, when written as a function of the configuration space and time, satisfies the Hamilton-Jacobi equation (see for instance~\cite{arnold}). The main goal of this section is showing that an analogous result holds for the stochastic action. 

Let $\varphi_t (\eta):\mathbb{D}_{t}\left(  \eta\right)  \rightarrow
M$ be the flow associated to the stochastic Hamilton equations (\ref{HJ 5}),
$(t,\eta)\in\mathbb{R}_{+}\times\Omega$. We define the function $R_{t}\left(
\eta\right)  :\mathbb{D}_{t}\left(  \eta\right)  \rightarrow\mathbb{R}$ as
$R_{t}\left(  \eta,z\right)  :=S\left(  \varphi\left(  z\right)  \right)
_{t}(\eta)$. The next proposition provides the differential of $R_{t}\left(  \eta\right)  $.

\begin{proposition}
\label{HJ th 1}Let $t\in\mathbb{R}_{+}$ be a fixed time instant and $\eta
\in\Omega$. Then $R_{t}\left(  \eta\right)  :\mathbb{D}_{t}\left(
\eta\right)  \rightarrow\mathbb{R}$ is differentiable and%
\begin{equation}
\label{derivative of r}
\mathbf{d}R_{t}\left(  \eta\right)  =\varphi_{t}\left(  \eta\right)  ^{\ast
}\theta-\theta,
\end{equation}
where $\theta$ is the one form of the exact symplectic manifold
$\left(  M,\omega=-\mathbf{d}\theta\right)  $.
\end{proposition}

\noindent\textbf{Proof.\ \ } We will proceed by showing that for any pair of points $x,y \in \mathbb{D}_{t}\left(  \eta\right)  $ we can write
\[
R_{t}\left(  \eta,x\right)  -R_{t}\left(  \eta,y\right)=\int_{\gamma}\left(  \varphi_{t}\left(  \eta\right)  ^{\ast}\left(
\theta\right)  -\theta\right),
\]
where $\gamma:\left(  a,b\right)  \subseteq
\mathbb{R}\rightarrow\mathbb{D}_{t}\left(  \eta\right)   $ is any smooth  curve in $\mathbb{D}_{t}\left(  \eta\right) $ that links $x $ and $y$.
This expression immediately implies that $R _t$ has continuous directional derivatives and it is hence Fr\'echet differentiable. Indeed, using first (\ref{HJ 6}), we have
\begin{equation}
\int_{\gamma}\left(  \varphi_{t}\left(  \eta\right)  ^{\ast}\left(
\theta\right)  -\theta\right)   =\int_{\gamma}\left(  \sum_{i=1}^{r}
\int_{0}^{t}\varphi_{s}^{\ast}\left(  \pounds _{X_{h_{i}}}\theta\right)
\delta X_{s}^{i}\right)  \left(  \eta\right)=\left(  \sum_{i=1}^{r}\int_{0}^{t}\left(  \int_{\gamma}\varphi_{s}^{\ast
}\left(  \pounds _{X_{h_{i}}}\theta\right)  \right)  \delta X_{s}^{i}\right)
\left(  \eta\right),  \label{HJ 7}
\end{equation}
where in the second equality we used Fubini's Theorem. Now, since
${\bf i}_{X_{h_{i}}}\omega= \mathbf{d}h_{i}$, for any $i=1,\ldots,r$, (\ref{HJ 7}) equals
\begin{align*}
&  \sum_{i=1}^{r}\left(  \int_{0}^{t}\left(  \int_{\gamma}\varphi_{s}^{\ast
}\mathbf{d}\left(  {\bf i}_{X_{h_{i}}}\theta\right)  \right)  \delta X_{s}^{i}%
-\int_{0}^{t}\left(  \int_{\gamma}\varphi_{s}^{\ast}\mathbf{d}h_{i}\right)
\delta X_{s}^{i}\right)  \left(  \eta\right) \\
&  =\sum_{i=1}^{r}\left(  \int_{0}^{t}\left(  \int_{\gamma}\mathbf{d}\left(
\varphi_{s}^{\ast}({\bf i}_{X_{h_{i}}}\theta)\right)  \right)  \delta X_{s}^{i}%
-\int_{0}^{t}\left(  \int_{\gamma}\mathbf{d(}\varphi_{s}^{\ast}h_{i})\right)
\delta X_{s}^{i}\right)  \left(  \eta\right) \\
&  =\sum_{i=1}^{r}\left(  \int_{0}^{t}\left[  {\bf i}_{X_{h_{i}}}\theta\left(
\varphi_{s}(\gamma_{b})\right)  -{\bf i}_{X_{h_{i}}}\theta\left(  \varphi_{s}\left(
\gamma_{a}\right)  \right)  \right]  \delta X_{s}^{i}-\int_{0}^{t}\left[
h_{i}\left(  \varphi_{s}(\gamma_{b})\right)  -h_{i}\left(  \varphi_{s}%
(\gamma_{a})\right)  \right]  \delta X_{s}^{i}\right)  \left(  \eta\right) \\
&  =\left(  \int_{0}^{t}\left\langle \theta,\delta\varphi_{s}\left(
\gamma_{b}\right)  \right\rangle -\int_{0}^{t}\left\langle \hat{h}\left(
\varphi_{s}\left(  \gamma_{b}\right)  \right)  ,\delta X_{s}\right\rangle
\right)  \left(  \eta\right)  -\left(  \int_{0}^{t}\left\langle \theta
,\delta\varphi_{s}\left(  \gamma_{a}\right)  \right\rangle -\int_{0}%
^{t}\left\langle \hat{h}\left(  \varphi_{s}\left(  \gamma_{a}\right)  \right)
,\delta X_{s}\right\rangle \right)  \left(  \eta\right) \\
&  =R_{t}\left(  \eta,x\right)  -R_{t}\left(  \eta,y\right).
\end{align*}
Given that $\gamma:\left(  a,b\right)  \rightarrow\mathbb{D}_{t}\left(
\eta\right)  $ and the points $x, y \in \mathbb{D}_{t}\left(
\eta\right)  $ are arbitrary, the result follows. $\ \ \ \ \blacksquare
$\bigskip

Later on in this section we will need the composition of $R$ 
with the inverse of the stochastic flow $\varphi$. More specifically, let
$\left(  t,\eta\right)  \in\mathbb{R}_{+}\times\Omega$ and let $\varphi
_{t}^{-1}\left(  \eta\right)  :\varphi_{t}\left(  \eta\right)  \left(
\mathbb{D}_{t}\left(  \eta\right)  \right)  \rightarrow\mathbb{D}_{t}\left(
\eta\right)  $ the inverse of $\varphi_{t}\left(  \eta\right)  $. We
define $\hat{R}_{t}\left(  \eta\right)  :\varphi_{t}\left(  \eta\right)
\left(  \mathbb{D}_{t}\left(  \eta\right)  \right)  \rightarrow\mathbb{D}
_{t}\left(  \eta\right)  $ as $\hat{R}_{t}\left(  \eta\right)  :=R_{t}\left(
\eta\right)  \circ\varphi_{t}^{-1}\left(  \eta\right)  =\varphi_{t}
^{-1}\left(  \eta\right)  ^{\ast}\left(  R_{t}\left(  \eta\right)  \right)  $.
Consequently,
\begin{equation}
\mathbf{d}\hat{R}_{t}\left(  \eta\right)  =\varphi_{t}^{-1}\left(
\eta\right)  ^{\ast}\left(  \mathbf{d}R_{t}\left(  \eta\right)  \right)
=\varphi_{t}^{-1}\left(  \eta\right)  ^{\ast}\left(  \varphi_{t}\left(
\eta\right)  ^{\ast}\left(  \theta\right)  -\theta\right)  =\theta-\varphi
_{t}^{-1}\left(  \eta\right)  ^{\ast}(\theta) \label{HJ 8}.
\end{equation}

In order to get closer to the classical deterministic result on the Hamilton-Jacobi equation we are first going to visualize it, using the map $R$, as a process depending on $M$ through the initial condition of the flow $\varphi$ generated by (\ref{HJ 5}). Second, we will restrict $R$ to a Lagrangian submanifold of $M$; this encodes mathematically the writing of the action as a function of the configuration space. Recall that a submanifold $\iota:L\hookrightarrow M$
of a symplectic manifold $(M,\omega)$ is called {\bfseries\itshape Lagrangian}
if $\dim\left(  L\right)  =\left.  \dim\left(  M\right)  \right/  2$ and
$\iota^{\ast}\omega=0$. Observe that since $\varphi_{t}\left(  \eta\right)  $
is a symplectomorphism a.s. for any $t\in\mathbb{R}_{+}$ and $\mathbb{D}
_{t}\left(  \eta\right)  $ is an open set, if $L$ is a Lagrangian submanifold
so are $L\cap\mathbb{D}_{t}\left(  \eta\right)  $ and $\varphi_{t}\left(
\eta\right)  \left(  L\cap\mathbb{D}_{t}\left(  \eta\right)  \right)  $. 

From now on we are going to assume that the underlying symplectic manifold $(M, \omega)$ is actually a cotangent bundle endowed with its canonical symplectic structure.
More specifically, $M=T^{\ast}Q$ for some manifold $Q$. In this case, a point $y\in
L\subset T^{\ast}Q$ in a Lagrangian submanifold $L$ is said to be a
{\bfseries\itshape regular point} of  $L$, if the restriction $\left.  \pi\right\vert
_{L}:L\rightarrow Q$ of the canonical projection $\pi:T^{\ast}Q\rightarrow Q$ to $L$  is a local diffeomorphism at $y$ (that is, $T_{y}\left.  \pi\right\vert
_{L}:T_{y}L\rightarrow T_{\pi(y)}Q$ is an isomorphism). In a
neighborhood $U\subset  L$ of a regular point $y\in L$ we can obviously describe the Lagrangian
submanifold $L$ using local coordinates on the base manifold $Q$,
which we will generally denote by $(  q^{1}, \ldots, q ^n)$.
On the other hand, since $\iota^{\ast}\omega=\mathbf{d}(\iota^{\ast}\theta
)=0$, there exists by the Poincar\'e lemma (shrinking $U$ if necessary) a smooth function $f\in C^{\infty}\left(  U\right)  $ such
that $\iota^{\ast}
\theta=\mathbf{d}f$. Conversely, if $\left(  q^{1}, \ldots, q ^n,p_{1}, \ldots, p _n
\right)  $ are local Darboux coordinates in a neighborhood
$V\subseteq T^{\ast}Q$ and $f\in C^{\infty}\left( \pi( V)\right)  $ is a function with no critical points, then the set
\begin{equation}
\label{lagrangian l f}
L_{f}=\left\{  \left(  q,p\right)  \in V\mid p_{i}=\frac{\partial f}{\partial
q^{i}},~i=1,...,n\right\}
\end{equation}
is a local Lagrangian submanifold such that 
\begin{equation}
\label{equation for the liouville form}
\iota _f ^\ast \theta= \pi |_{L_f} ^\ast \mathbf{d}f, 
\end{equation}
with $\iota_f: L _f \hookrightarrow V $ the inclusion and $ \pi |_{L_f}: L _f\subset T ^\ast Q \rightarrow \pi(V)  $ the local diffeomorphism obtained by restriction of the canonical projection. 

\begin{theorem}
\label{theorem 1}
Let $Q$ be a manifold and let $L\subset
T^{\ast}Q$ a Lagrangian submanifold. Let $y_{0}\in L$ be a regular point and
let $x_{0}=\pi(y_{0})$, where $\pi:T^{\ast}Q\rightarrow Q$ is the canonical
projection. Then, there exist two neighborhoods $V_{y_{0}}\subseteq L$ and $V_{x_{0}}\subseteq
Q$  of $y_{0}$ and of $x_{0}$, respectively and a
map $\xi:\Omega\times V_{x_{0}}\rightarrow\mathbb{R}_{+}$ with the property
that $\xi\left(  x\right)  :\Omega\rightarrow\mathbb{R}_{+}$ is a stopping
time, such that the equation%
\begin{equation}
\pi\left(  \varphi_{s}\left(  \eta,y\right)  \right)  =x\label{HJ 12}%
\end{equation}
has a unique solution in $V_{y_{0}}\subseteq L$ for any $\eta\in\Omega$, any
$x\in V_{x_{0}}$, and any $s\in\lbrack0,\xi\left(  \eta,x\right)  ]$. We are
going to denote this solution by $\psi_{s}\left(  \eta,x\right)  $. Moreover,
$\psi\left(  x\right)  :[0,\xi(x))\rightarrow V_{y_{0}}\left(  \eta\right)  $
is a semimartingale for any $x\in V_{x_{0}}$ and $\psi_{s}\left(  \eta\right)
:V_{x_{0}}\rightarrow V_{y_{0}}$ is a diffeomorphism for any $s\in\lbrack
0,\xi\left(  x\right)  )$ which depends continuously on $s$.
\end{theorem}

\noindent\textbf{Proof.} \ \ Let $U_{y_{0}}\subset L$ be an open neighborhood of
$y_{0}\in L$. We pick $U_{y _0}$ small enough so that $\pi| _{U_{y_{0}}} $ is a diffeomorphism onto its image and a set of local coordinates $\left(  q^{i};i=1,...,n\right)  $ can be chosen 
on $U_{x_{0}}:=\pi\left(  U_{y_{0}}\right)  $. Let $\left(  y^{i}
=q^{i}\circ\left.  \pi\right\vert _{L};i=1,...,n\right)  $ be the corresponding
induced coordinates on $U_{y_{0}}$. Denote by $\hat{q}:U_{x_{0}}
\rightarrow\mathbb{R}^{n}$ and $\hat{y}:U_{y_{0}}\rightarrow\mathbb{R}^{n}$
the local chart maps associated to these coordinates. For any $y\in U_{y_{0}}$,
let $\tau_{U_{x_{0}}}\left(  y,\eta\right)  =\inf\{t>0~|~\pi\circ\varphi
_{t}\left(  \eta,y\right)  \notin U_{x_{0}}\}$ be the first exit time at which the
semimartingale $\pi\circ\varphi\left(  y\right)  $ leaves $U_{x_{0}}$. Let $F$
be the restriction of $\pi\circ\varphi$ to the set $A:=\{\left(
s,\eta,y\right)  \in\mathbb{R}_{+}\times\Omega\times U_{y_{0}}~|~s\in
\lbrack0,\tau_{U_{x_{0}}}\left(  y,\eta\right)  )\}$. In local coordinates,
$F:A\rightarrow U_{x_{0}}$ is expressed as%
\[
F_{s}^{j}\left(  \eta\right)  \left(  y^{1},...,y^{n}\right)  =q^{j}\circ
\pi\circ\varphi_{s}\left(  \eta\right)  \circ\hat{y}^{-1}\left(
y^{1},...,y^{n}\right)  \text{, \ }j=1,...,n.
\]
Now, remark that $\det\left(  \frac{\partial F_{0}^{j}(\eta)}{\partial y^{i}%
}\left(  y_{0}\right)  \right)  \neq0$ a.s. because $y_{0}\in L$ is a regular
point. The continuity of the derivative of $F_{0}\left(  \eta\right)  :U_{y_{0}}\rightarrow U_{x_{0}}$ implies that there exists a neighborhood $V_{y_{0}}\subseteq U_{y_{0}}$ such
that $\det\left(  \frac{\partial F_{0}^{j}(\eta)}{\partial y^{i}}\left(
y\right)  \right)  >0$ a.s., for any $y\in V_{y_{0}}$. For any of these $y\in
V_{y_{0}}$, let%
\[%
\begin{array}
[c]{rrl}%
Z\left(  y\right)  :=\det\left(  \frac{\partial F^{j}}{\partial y^{i}}\left(
y\right)  \right)  :\left[  0,\tau_{U_{x_{0}}}\left(  y\right)  \right)   &
\longrightarrow & \mathbb{R}\\
\left(  s,\eta\right)   & \longmapsto & \det\left(  \frac{\partial F_{s}%
^{j}(\eta)}{\partial y^{i}}\left(  y\right)  \right),
\end{array}
\]
which is a well defined and continuous semimartingale, by the continuity of the
differential of the flow $\varphi$. Observe that $Z_{0}\left(  y\right)  >0$ for any $y\in V_{y_{0}}$.
Let $T\left(  y,\eta\right)  :=\inf\{\tau_{U_{x_{0}}}\left(  y\right)  \geq
t>0~|~Z_{t}\left(  y,\eta\right)  \notin\mathbb{R}_{+}\}$.

Now, recall that we want to see that the equation $\pi\left(  \varphi_{s}\left(
\eta,y\right)  \right)  =x$ has a unique solution in $y\in L$, for any $x\in
V_{x_{0}}$ in a suitable $V_{x_{0}}$  and up to a suitable stopping time $\xi\left(  x\right)  $. Therefore, it suffices to solve the equation
\begin{equation}
\pi\left(  \varphi_{s}^{T(y)}\left(  \eta,y\right)  \right)  =x,\label{HJ 9}%
\end{equation}
where $\varphi^{T\left(  y\right)  }\left(  y\right)  $ denotes
the process $\varphi\left(  y\right)  $ stopped at  time $T\left(
y\right)  $, that is, $\varphi^{T\left(  y\right)  }\left(  y\right)  \left(
s,\eta\right)  =\varphi_{T\left(  y,\eta\right)  \wedge s}\left(
\eta,y\right)  $. Observe that $\varphi^{T\left(  y\right)  }\left(  y\right)
$ is always in $U_{y_{0}}$ if $y$ was already in $V_{y_{0}}$. Consequently,
$\varphi^{T\left(  y\right)  }\left(  y\right)  $ may be described using the local
coordinates introduced above. Moreover, if we set $\xi\left(  x\right) : =T\left(  \left.  \pi\right\vert_{L}^{-1}(x)\right)  $, $V_{x_{0}}:=\pi(V_{y_{0}})$, the equation~(\ref{HJ 9}) admits by construction a unique solution $\psi_{s}\left(  \eta,x\right)  $ via the Implicit Function Theorem.
Additionally, if we apply
the Stratonovich differentiation rules to%
\[
\pi\left(  \varphi_{s}^{T(y)}\left(  \eta,\psi_{s}\left(  \eta,x\right)
\right)  \right)  =x\text{, \ }s\in\lbrack0,\xi\left(  x,\eta\right)  )
\]
we obtain that $\psi_{s}(\eta,x)$ satisfies up to time $\xi\left(  x\right)  $ the  Stratonovich
differential equation
\begin{equation}
\delta\psi_{s}(x)= \sum
_{i=1}^{r}\left[  T_{\psi_{s}(x))}F\right]  ^{-1}\left(
T_{\varphi^{T\left(  y\right)
}_{s}\left(  \psi_{s}(x)\right)  }(\hat{q}\circ\pi)\left( X_{h_{i}}(\varphi^{T\left(  y\right)
}_{s}\left(  \psi_{s}(x)\right)  )\right)  \right)\delta
X_{s}^{i}  \label{HJ 11}%
\end{equation}
with initial condition $\psi_{s=0}(x)=y(x)\in V_{y_{0}}$ a.s. such that
$\pi\left(  y(x)\right)  =x\in V_{x_{0}}$. That is, we can visualize $\psi_{s}(\eta,x)$ as the unique stochastic flow associated to the stochastic
differential equation (\ref{HJ 11}). This guarantees that the properties claimed in the statement hold. \ \ \ \ $\blacksquare$

\medskip

We proceed now by considering  the stochastic action $R$ not
as a semimartingale parametrized by $T^{\ast}Q$ through the initial condition
of the stochastic flow $\varphi$ defined by (\ref{HJ 5}), but as a process
depending on the base manifold $Q$. More specifically, we will restrict to the open neighborhood $V_{x_{0}} \subset  Q$ introduced in the statement of
Theorem \ref{theorem 1} and which is
mapped onto $V_{y_{0}} \subset L$ using the map $\psi$ that solves (\ref{HJ 12}).
Furthermore, since we are always going to work around regular points of the Lagrangian submanifold, we will always consider  Lagrangian
submanifolds of the type $L_{f}$  (see~(\ref{lagrangian l f})) for some $f\in C^{\infty}\left(  Q\right)  $.

\begin{definition}
Let $L_{f}\subseteq T^{\ast}Q$ be a Lagrangian submanifold, $f\in C^{\infty
}(Q)$. Let $V_{x_{0}}\subseteq Q$ be the open neighborhood of $x_{0}$ introduced
in Theorem \ref{theorem 1} and $\psi(x):[0,\xi(x))\rightarrow V_{y_{0}}$ the
semimartingale solution of (\ref{HJ 11}) with initial condition $x\in
V_{x_{0}}$ a.s.. We define the {\bfi  projected stochastic action} $\widetilde{S}\left(  x\right)  :[0,\xi
(x))\rightarrow\mathbb{R}$ as%
\[
\widetilde{S}_{t}\left(  \eta,x\right)  =R_{t}\left(  \eta,\psi_{t}%
(\eta,x)\right)  +f\left(  \pi\left(  \psi_{t}(\eta,x)\right)  \right)
=\left(  R_{t}\left(  \eta\right)  +f\circ\pi\right)  \circ\psi_{t}(\eta,x).
\]

\end{definition}

Notice that the differentiability properties of the maps
$R$, $f\in C^{\infty}\left(  Q\right)  $, and  $\psi$ imply that 
the map
\begin{equation}
\begin{array}[c]{rrl}
\widetilde{S}_{t}\left(  \eta\right)  :\mathbb{D}_{t}^{\psi}(\eta) &
\longrightarrow & \mathbb{R}\\
x & \longmapsto & \widetilde{S}_{t}\left(  \omega,x\right)
\end{array}
\label{HJ 13}
\end{equation}
is continuously differentiable for any $\left(  t,\eta\right)  \in\mathbb{R}
_{+}\times\Omega$ such that $t\in\lbrack0,\xi\left(  x,\eta\right)  )$. In this expression $\mathbb{D}_{t}^{\psi
}\left(  \eta\right) : =\{x\in V_{x_{0}}~|~t<\xi\left(  x,\eta\right)  \}$. 
The following theorem provides an explicit expression for the spatial derivatives
of the projected stochastic action $\widetilde{S}$.

\begin{theorem}
\label{HJ th 3}
Let  $L _f $ be a Lagrangian submanifold of $T ^\ast Q $,  $f \in C^\infty(Q)$. Then, on  the open set $\mathbb{D}_{t}^{\psi}(\eta)$,  $\left(  t,\eta\right)  \in\mathbb{R}_{+}\times\Omega$,
\begin{equation}
\label{derivative s tilde}
\mathbf{d}\widetilde{S}_{t}\left(  \eta\right)  =\left(  \varphi_{t}\left(
\eta\right)  \circ\psi_{t}\left(  \eta\right)  \right)  ^{\ast}\theta.
\end{equation}
If $\left(
q^{i},p_{i};i=1,...,n\right)  $ are local Darboux coordinates of $T^{\ast}Q$ on
an open neighborhood of a regular point $y_{0}\in L_{f}$, the expression~(\ref{derivative s tilde}) can be locally written as
\[
\frac{\partial\widetilde{S}_{t}(\eta)}{\partial q^{i}}\left(  q\right)
=p_{i}\left(  \varphi_{t}\left(  \eta,\psi_{t}(\eta,q)\right)  \right)
,~~i=1,...,n.
\]

\end{theorem}

\noindent\textbf{Proof.} \ \ First of all observe that $\widetilde{S}%
_{t}\left(  \eta\right)  $ can be expressed in terms of $\hat{R}_{t}\left(
\eta\right)  $ as follows:%
\[
\widetilde{S}_{t}\left(  \eta,q\right)  =\hat{R}_{t}\left(  \eta\right)
\circ\varphi_{t}\left(  \eta\right)  \circ\psi_{t}\left(  \eta,q\right)
+f\circ\pi\circ\psi_{t}\left(  \eta,q\right)  .
\]
Then, for any smooth curve $\gamma:[a,b] \rightarrow \mathbb{D} _t ^\psi (\eta) $
\begin{equation}
\widetilde{S}_{t}\left(  \eta,\gamma_{b}\right)  -\widetilde{S}_{t}\left(
\eta,\gamma_{a}\right)  =\int_{\gamma}\mathbf{d}\widetilde{S}_{t}\left(
\eta\right)  =\int_{\gamma}\mathbf{d}\left[  \hat{R}_{t}\left(  \eta\right)
\circ\varphi_{t}\left(  \eta\right)  \circ\psi_{t}\left(  \eta\right)
\right]  +\int_{\gamma}\mathbf{d}\left(  f\circ\pi\circ\psi_{t}\left(
\eta\right)  \right)  .\label{HJ 14}%
\end{equation}
Given that $\mathbb{D} _t ^\psi (\eta) \subset L _f  $, the curve $\gamma  $ takes values in the Lagrangian submanifold $L _f $ and hence~(\ref{HJ 14}) can be rewritten as
\begin{eqnarray}
\label{HJ 14 bis}
\widetilde{S}_{t}\left(  \eta,\gamma_{b}\right)  -\widetilde{S}_{t}\left(
\eta,\gamma_{a}\right)&=&\int_{\gamma} \iota ^\ast _{L _f}\mathbf{d}\left[  \hat{R}_{t}\left(  \eta\right)
\circ\varphi_{t}\left(  \eta\right)  \circ\psi_{t}\left(  \eta\right)
\right]  +\int_{\gamma}\mathbf{d}\left(  f\circ\pi\circ\psi_{t}\left(
\eta\right)  \right)  \notag\\
	&= &\int_{\gamma}\mathbf{d}\left[  \hat{R}_{t}\left(  \eta\right)
\circ\varphi_{t}\left(  \eta\right)  \circ\psi_{t}\left(  \eta\right)\circ  \iota _{L _f}
\right]  +\int_{\gamma}\mathbf{d}\left(  f\circ\pi\circ\psi_{t}\left(
\eta\right)  \right). \label{almost hj8}
\end{eqnarray}
On the other hand, we saw in (\ref{HJ 8}) that%
\[
\mathbf{d}\hat{R}_{t}=\theta-\varphi_{t}^{-1}\left(  \eta\right)  ^{\ast
}(\theta).
\]
Moreover, since $\iota _f ^\ast \theta= \pi |_{L_f} ^\ast \mathbf{d}f $ we have that
\begin{equation*}
\mathbf{d}\left[  \hat{R}_{t}\left(  \eta\right)
\circ\varphi_{t}\left(  \eta\right)  \circ\psi_{t}\left(  \eta\right)\circ  \iota _{L _f}
\right]= \left(\varphi_{t}\left(  \eta\right)  \circ\psi_{t}\left(  \eta\right)\circ  \iota _{L _f} \right) ^\ast \theta- \mathbf{d} \left( f \circ \pi \circ \psi _t(\eta)\right),
\end{equation*}
which substituted in~(\ref{almost hj8}) yields
\[
\widetilde{S}_{t}\left(  \eta,\gamma_{b}\right)  -\widetilde{S}_{t}\left(
\eta,\gamma_{a}\right)=\int _\gamma \left(\varphi_{t}\left(  \eta\right)  \circ\psi_{t}\left(  \eta\right)\circ  \iota _{L _f} \right) ^\ast \theta=\int _\gamma \left(\varphi_{t}\left(  \eta\right)  \circ\psi_{t}\left(  \eta\right) \right) ^\ast \theta.
\]
Since $\gamma$ is an arbitrary smooth curve, we can conclude that
\[
\mathbf{d}\widetilde{S}_{t}\left(  \eta\right)  =\left(  \varphi_{t}\left(
\eta\right)  \circ\psi_{t}\left(  \eta\right)  \right)  ^{\ast}\theta,
\]
as required. \ \ \ \ $\blacksquare$

\medskip

We conclude this section by proving that the projected stochastic action $\widetilde{S}_{t}$ satisfies a specific stochastic differential equation which generalizes the classical Hamilton-Jacobi equation. For obvious reasons, this equation will be referred to as the
{\bfi stochastic Hamilton-Jacobi equation}.

\begin{theorem}[Stochastic Hamilton-Jacobi equation]
\label{theorem Hamilton-Jacobi}
Using the same notation as in Theorem
\ref{theorem 1}, the projected stochastic action $\widetilde{S}\left(  q\right)
:[0,\xi\left(  q\right)  )\rightarrow\mathbb{R}$ associated to the Lagrangian submanifold $L _f  $  defined by the function $f \in C^\infty(Q) $ satisfies
\[
\widetilde{S}\left(  q\right)  =f\left(  q\right)  -\int\left\langle\hat{h}\left(
q,\frac{\partial\widetilde{S}_{s}}{\partial q}\left(  q\right)  \right)
,\delta X_{s}\right\rangle
\]
for any $q\in V_{x_{0}}$.
\end{theorem}

In order to prove this theorem we need the following auxiliary result.

\begin{proposition}
[{\cite[Theorem 3.3.2]{Kunita book}}]\label{HJ prop 2} Let $F(x):\mathbb{R}%
_{+}\times\Omega\rightarrow\mathbb{R}$, $x\in\mathbb{R}^{n}$, be a family of continuous
semimartingales parametrized by $\mathbb{R} ^n $. Suppose that the dependence of this family on the $\mathbb{R}^n $parameter  is at least three times differentiable. In addition, suppose that there exists a process $f:\mathbb{R}_{+}\times\Omega\times\mathbb{R}%
^{n}\rightarrow\mathbb{R}^{d}$ that satisfies sufficient regularity conditions and
a semimartingale $X:\mathbb{R}_{+}\times\Omega\rightarrow\mathbb{R}^{d}$
such that
\[
F(x)=\sum_{j=1}^{r}\int f_{j}\left(  t,x\right)  \delta X_{t}^{j}.
\]
Let $g:\mathbb{R}_{+}\times\Omega\rightarrow\mathbb{R}^{n}$ be a continuous 
$\mathbb{R}^{n}$-valued semimartingale. Then $F\left(  g\right)
:\mathbb{R}_{+}\times\Omega\rightarrow\mathbb{R}$ defined as $F\left(
g\right)  \left(  t,\eta\right)  :=F\left(  g_{t}\left(  \eta\right)
,t,\eta\right)  $ satisfies%
\[
F\left(  g_{t},t\right)  -F\left(  g_{0},0\right)  =\sum_{j=1}^{r}\int_{0}%
^{t}f_{j}\left(  s,g_{s}\right)  \delta X_{s}^{j}+\sum_{i=1}^{n}\int_{0}%
^{t}\frac{\partial F}{\partial x^{i}}\left(  s,g_{s}\right)  \delta g_{s}%
^{i}.
\]

\end{proposition}

\noindent\textbf{Proof of Theorem \ref{theorem Hamilton-Jacobi}.} \ \ First of
all observe that using the definition of the function $R_{t}$ the semimartingale $\widetilde{S}\left(  q\right)
:[0,\xi\left(  q\right)  )\rightarrow\mathbb{R}$ may be expressed as
\begin{align*}
\widetilde{S}\left(  q\right)   &  =f\circ\pi\circ\psi_{t}\left(
\eta,q\right)  +R_{t}\left(  \eta,\psi_{t}(\eta,q)\right)  \\
&  =f\circ\pi\circ\psi_{t}\left(  \eta,q\right)  +\left.  \sum_{j=1}%
^{r}\left(  \int\left( {\bf  i}_{X_{h_{j}}}\theta-h_{j}\right)  \left(  \varphi
_{s}(z)\right)  \delta X_{s}^{j}\right)  \right\vert _{z=\psi_{t}(\eta,q)}.
\end{align*}
If we use Proposition~\ref{HJ prop 2} in the second summand of this expression, we obtain
\begin{equation}
\widetilde{S}\left(  q\right)  =f\circ\pi\circ\psi\left(  q\right)
+\sum_{j=1}^{r}\left(  \int\left(  {\bf i}_{X_{h_{j}}}\theta-h_{j}\right)  \left(
\varphi_{s}(\psi_{s}(q))\right)  \delta X_{s}^{j}\right)  +\int\left\langle
\mathbf{d}R_{s},\delta\psi_{s}\left(  q\right)  \right\rangle .\label{HJ 10}%
\end{equation}
We now separately study  the summands in the right hand side
of this equation in order to prove the statement of the theorem. We start by recalling that by
Proposition \ref{HJ th 1}, $\mathbf{d}R_{s}=\varphi_{s}^{\ast}\theta-\theta$ and hence
\begin{equation}
\int\left\langle \mathbf{d}R_{s},\delta\psi_{s}\left(  q\right)  \right\rangle
=\int\left\langle \varphi_{s}^{\ast}\theta-\theta,\delta\psi_{s}\left(
q\right)  \right\rangle \text{.}\label{HJ 30}%
\end{equation}
Furthermore, since $\iota _f ^\ast \theta= \pi |_{L_f} ^\ast \mathbf{d}f$ and the
semimartingale $\psi\left(  q\right)  $ takes values in $V_{y_{0}}\subseteq
L_f$,%
\begin{equation}
\int_{0}^{t}\left\langle \theta,\delta\psi_{s}\left(  q\right)  \right\rangle
=\int_{0}^{t}\left\langle \mathbf{d}\left(  f\circ\pi\right)  ,\delta\psi
_{s}\left(  q\right)  \right\rangle =f\circ\pi\circ\psi_{t}\left(  q\right)
-f\left(  q\right)  .\label{HJ 15}%
\end{equation}
We now recall that the semimartingale $\varphi(\psi(q)):[0,\xi(q))\rightarrow
T^{\ast}Q$ takes values in the fiber $\pi^{-1}\left(  q\right)  $. Indeed,
by the construction in Theorem~\ref{theorem 1}, $\psi\left(  q\right)  $ is the semimartingale starting at $q$
such that%
\[
\pi\left(  \varphi_{s}\left(  \eta,\psi_{s}\left(  \eta,q\right)  \right)
\right)  =q
\]
for any $\left(  s,\eta\right)  \in\lbrack0,\xi(q))$. Then, since $\theta$ is
a semibasic form we necessarily have that
\[
\int\left\langle \theta,\delta\left(  \varphi_{s}(\psi_{s}\left(  q\right)
)\right)  \right\rangle =0.
\]
But, using the fact that $\varphi$ is the flow of the stochastic Hamilton
equations (\ref{HJ 5}), by Proposition \ref{HJ prop 2}, we have that for any
$g\in C^{\infty}\left(  M\right)  $%
\begin{equation}
g\left(  \varphi\left(  \psi\left(  q\right)  \right)  \right)  =g\left(
y\left(  q\right)  \right)  +\sum_{j=1}^{r}\int X_{h_{j}}[g](\varphi
_{s}\left(  \psi_{s}\left(  q\right)  \right)  )\delta X_{s}^{j}%
+\int\left\langle \mathbf{d}\left(  g\circ\varphi_{s}\right)  ,\delta\psi
_{s}(q)\right\rangle \label{HJ 16}%
\end{equation}
where $y\left(  q\right)  \in L_{f}$ is the unique point such that $\left.
\pi\right\vert _{L}(y(q))=q$. We claim that
\begin{equation}
0=\int\left\langle \theta,\delta\left(  \varphi_{s}(\psi_{s}\left(  q\right)
)\right)  \right\rangle =\sum_{j=1}^{r}\int\left( {\bf  i}_{X_{h_{j}}}\theta\right)
\left(  \varphi_{s}\left(  \psi_{s}\left(  q\right)  \right)  \right)  \delta
X_{s}^{j}+\int\left\langle \varphi_{s}^{\ast}\theta,\delta\psi_{s}%
(q)\right\rangle .\label{HJ 17}%
\end{equation}
Indeed, since we are working at a local level we can use Darboux coordinates and we can replace $\theta$ by $\sum_{i=1}^{n}p_{i}\mathbf{d} q^{i}$; (\ref{HJ 17}) is a straightforward
consequence of (\ref{HJ 16}). If we now plug (\ref{HJ 30}), (\ref{HJ 15}), and
(\ref{HJ 17}) into (\ref{HJ 10}) we obtain
\begin{equation}
\widetilde{S}\left(  q\right)  =f\left(  q\right)  -\sum_{j=1}^{r}\int
h_{j}\left(  \varphi_{s}\left(  \psi_{s}(q)\right)  \right)  \delta X_{s}
^{j}.\label{HJ 18}%
\end{equation}
Finally, we saw in Theorem \ref{HJ th 3} that
\[
p_{i}\left(  \varphi_{t}\circ\psi_{t}\left(  \eta,q\right)  \right)
=\frac{\partial\widetilde{S}_{t}\left(  \eta\right)  }{\partial q^{i}}\left(
q\right)  ,~~i=1,...,n,
\]
on $\mathbb{D}_{t}^{\psi}\left(  \eta\right)  =\{x\in V_{x_{0}}~|~\xi\left(
x,\eta\right)  >t\}$, $\left(  t,\eta\right)  \in\mathbb{R}_{+}\times\Omega$.
For any $\eta\in\Omega$, the time parameter $s$ in the integrand of
(\ref{HJ 18}) is always smaller than $\xi\left(  q,\eta\right)  $ and hence as
$\frac{\partial \widetilde S_s}{\partial q^{i}}\left(  q\right)  $ and
$p_{i}\left(  \varphi_{s}\circ\psi_{s}\left(  q\right)  \right)  $ coincide
a.s. on $[0,\xi(q))$ for any $i=1,...,n$, the result follows.
\ \ \ \ $\blacksquare$

\begin{example}
\label{feynman kac action} \normalfont Let $Q=\mathbb{R}^{n}$ and $T^{\ast
}Q=\mathbb{R}^{n}\times\mathbb{R}^{n}$ with global coordinates $\left(
q^{i},p_{i};i=1,...,n\right)  $. Let $f\in C^{\infty}\left(  \mathbb{R}%
^{n}\right)  $, $h_{0}\in C^{\infty}\left(  \mathbb{R}^{2n}\right)  $, and
$h_{i}=p_{i}$ for any $i=1,...,n$. Consider the semimartingale $X:\mathbb{R}%
_{+}\times\Omega\rightarrow\mathbb{R}^{n+1}$ given by $\left(  t,\omega
\right)  \mapsto(t,B_{t}^{1},...,B_{t}^{n})$, where $(B^{1},...,B^{n})$ is an
$n$-dimensional Brownian motion. That is, $[B_{t}^{i},B_{t}^{j}]=\delta^{ij}%
t$, where $[\cdot,\cdot]$ denotes the quadratic variation. Then, the projected
stochastic action $\widetilde{S}:\mathbb{R}_{+}\times\Omega\times
\mathbb{R}^{n}\rightarrow\mathbb{R}$ built from the stochastic Hamiltonian
system on $\mathbb{R}^{2n}$ with Hamiltonian fuction $h=(h_{0},h_{1}%
,...,h_{n})$ and stochastic component $X$ satisfies by Theorem
\ref{theorem Hamilton-Jacobi}%
\begin{equation}
\widetilde{S}_{t}\left(  q\right)  =f(q)-\int_{0}^{t}h_{0}\left(
q,\frac{\partial\widetilde{S}_{s}}{\partial q}(q)\right)  ds-\sum_{i=1}%
^{n}\int\frac{\partial\widetilde{S}_{s}}{\partial q^{i}}(q)\delta B_{s}%
^{i}.\label{HJ 31}%
\end{equation}
If we transform the It\^{o} integrals in this expression into Stratonovich
integrals, (\ref{HJ 31}) reads%
\begin{equation}
\widetilde{S}_{t}\left(  q\right)  =f(q)-\int_{0}^{t}h_{0}\left(
q,\frac{\partial\widetilde{S}_{s}}{\partial q}(q)\right)  ds-\sum_{i=1}%
^{n}\int_{0}^{t}\frac{\partial\widetilde{S}_{s}}{\partial q^{i}}(q)dB_{s}%
^{i}-\frac{1}{2}\sum_{i=1}^{n}\left[  \frac{\partial\widetilde{S}}{\partial
q^{i}}(q),B^{i}\right]  _{t},\label{HJ 39}%
\end{equation}
It is not difficult to realize, though tedious to check, that%
\[
\frac{\partial\widetilde{S}_{t}}{\partial q^{i}}\left(  q\right)
=\frac{\partial f}{\partial q^{i}}(q)-\int_{0}^{t}\frac{\partial}{\partial
q^{i}}\left(  h_{0}\left(  q,\frac{\partial\widetilde{S}_{s}}{\partial
q}(q)\right)  \right)  ds-\sum_{r=1}^{n}\int_{0}^{t}\frac{\partial}{\partial
q^{i}}\left(  h_{r}\left(  q,\frac{\partial\widetilde{S}_{s}}{\partial
q}(q)\right)  \right)  \delta B_{s}^{r}%
\]
so that, since $h_{r}=p_{r}$ for any $r=1,...,r$,
\[
\frac{\partial\widetilde{S}_{t}}{\partial q^{i}}\left(  q\right)
=\frac{\partial f}{\partial q^{i}}(q)-\int_{0}^{t}\frac{\partial}{\partial
q^{i}}\left(  h_{0}\left(  q,\frac{\partial\widetilde{S}_{s}}{\partial
q}(q)\right)  \right)  ds-\sum_{r=1}^{n}\int_{0}^{t}\frac{\partial
^{2}\widetilde{S}_{s}}{\partial q^{i}\partial q^{r}}(q)\delta B_{s}^{r}.
\]
Therefore, disregarding all the finite variation terms in this last
expression, we have%
\begin{align*}
\left[  \frac{\partial\widetilde{S}}{\partial q^{i}}(q),B^{i}\right]  _{t}  &
=-\sum_{r=1}^{n}\left[  \int\frac{\partial^{2}\widetilde{S}_{s}}{\partial
q^{i}\partial q^{r}}(q)dB_{s}^{r},\int dB_{s}^{i}\right]  _{t}=-\sum_{r=1}%
^{n}\int_{0}^{t}\frac{\partial^{2}\widetilde{S}_{s}}{\partial q^{i}\partial
q^{r}}(q)d[B^{r},B^{i}]_{s}\\
& =-\sum_{r=1}^{n}\int_{0}^{t}\frac{\partial^{2}\widetilde{S}_{s}}{\partial
q^{i}\partial q^{r}}(q)\delta^{ir}ds=-\int_{0}^{t}\frac{\partial^{2}%
\widetilde{S}_{s}}{(\partial q^{i})^{2}}(q)ds,
\end{align*}
where the property%
\[
\left[  \int HdX,\int KdY\right]  _{t}=\int_{0}^{t}H_{s}K_{s}d[X,Y]_{s}%
\]
for arbitrary real semimartingales $H$, $K$, $X$, and $Y$ (\cite[Chapter II
Theorem 29]{protter book}) has been used. Taking expectations in both sides of
(\ref{HJ 39}) and assuming that all the processes involved are regular enough
so that Fubini's Theorem may be invoked, we obtain%
\[
E[\widetilde{S}_{t}\left(  q\right)  ]=f(q)-\int_{0}^{t}E\left[  h_{0}\left(
q,\frac{\partial\widetilde{S}_{s}}{\partial q}(q)\right)  \right]  ds+\frac
{1}{2}\int_{0}^{t}\Delta E[\widetilde{S}_{s}(q)]ds
\]
Finally, take $h_{0}=\frac{1}{2}\sum_{i=1}^{n}p_{i}^{2}+V(q)$, $V\in
C^{\infty}\left(  \mathbb{R}^{n}\right)  $, and let $\Phi_{t}(q):=\exp
(-E[\widetilde{S}_{t}\left(  q\right)  ])$. Then%
\begin{align*}
\frac{\partial}{\partial t}\Phi_{t}(q) &  =\Phi_{t}(q)\left[  V(q)+\frac{1}%
{2}\sum_{i=1}^{n}E\left[  \left(  \frac{\partial\widetilde{S}_{s}}{\partial
q^{i}}(q)\right)  ^{2}\right]  -\frac{1}{2}\Delta E[\widetilde{S}%
_{s}(q)]\right]  \\
&  =V(q)\Phi_{t}(q)+\frac{1}{2}\Delta\Phi_{t}(q).
\end{align*}
This shows that the projected stochastic action $\widetilde{S}_{t}$ can be
used to construct solutions of the heat equation modified with a potential
term $V$, with initial condition given by the function $\exp(f)\in C^{\infty
}\left(  \mathbb{R}^{n}\right)  $.
\end{example}

\section{The Hamilton-Jacobi equation and generating functions}

One of the main features of  the Hamilton-Jacobi equation is that its solutions can be used as generating functions of time-dependent symplectomorphisms that transform the original Hamiltonian system in such a way that its solutions can be easily written down. The goal of this section is spelling out the way in which this classical procedure can be extended to stochastic Hamiltonian systems. 

\subsection{The deterministic case}
\label{seccion HJ determinista}

We start by recalling the relation between the Hamilton-Jacobi equation and the generating functions for integrating canonical transformations in the classical deterministic case. In the next paragraphs we will write down some classical results in a form that is well adapted for the subsequent generalization to the stochastic case. All along this section we will consider Hamiltonian systems on cotangent bundles $\left(  T^{\ast}Q,\omega= - \mathbf{d} \theta\right)  $ endowed with their canonical symplectic forms.

Consider the manifold $T^{\ast}Q\times T^{\ast}Q$ endowed with the symplectic form $\Omega:=\tau_{1}^{\ast}\omega-\tau_{2}^{\ast}\omega$, where $\tau_{i}:T^{\ast}Q\times T^{\ast
}Q\rightarrow T^{\ast}Q$, $i=1,2$, denote the canonical projections onto the
first and the second factors, respectively. Let now $\psi:T^{\ast}Q\rightarrow T^{\ast}Q$ be a smooth function. It is easy to verify that the map $\psi$ is a
symplectomorphism if and only if $\iota_{\psi}^{\ast}\Omega=0$, where
$\iota_{\psi}:L^{\psi}\hookrightarrow T^{\ast}Q\times T^{\ast}Q$ is the
inclusion of the graph $L^{\psi}$ of $\psi$ (\cite[Proposition 5.2.1]{fom}),
in which case is a Lagrangian submanifold of $T^{\ast}Q\times T^{\ast}Q $. Given that
$\Omega=-\mathbf{d}\Theta$, with $\Theta=\tau_{1}^{\ast}\theta
-\tau_{2}^{\ast}\theta$, we have that $0=\iota_{\psi}^{\ast}\Omega
=-\iota_{\psi}^{\ast}\left(  \mathbf{d}\Theta\right)  =-\mathbf{d(}\iota
_{\psi}^{\ast}\Theta)$ and hence by Poincar\'{e}'s Lemma, we can locally write $\iota_{\psi}^{\ast}\Theta=\mathbf{d}S$, for some function $S\in C^{\infty}\left(
L^{\psi}\right)  $. We will say that  $S$ is a local {\bfseries\itshape generating
function} for the symplectic map $\psi$. In addition, suppose that
\begin{equation}
\tau:T^{\ast}Q\times T^{\ast}Q\rightarrow Q\times Q,~~\tau=\pi\circ\tau
_{1}\times\pi\circ\tau_{2}\label{HJ 23}%
\end{equation}
with $\pi:T^{\ast}Q\rightarrow Q$ the canonical projection, is a local
diffeomorphism when restricted to $L^{\psi}$ and denote its (local) inverse by
$\tau^{-1}:Q\times Q\rightarrow L^{\psi}$. We will suppose throughout this section that this is the case and we will think of  the generating function $S\in C^{\infty}\left(  L^{\psi
}\right)  $ as a function defined on $Q\times Q$; that is, we will not
distinguish between $S$ and $\left(  \tau^{-1}\right)  ^{\ast}S$. With this convention, we can write
\begin{equation}
\mathbf{d}_{Q\times Q}S=\left(  \tau^{-1}\right)  ^{\ast}\circ\iota_{\psi
}^{\ast}\left(  \Theta\right).  \label{HJ 24}%
\end{equation}

Let now $\{\psi_{t}\}_{t\in\mathbb{R}}$ be a family of symplectomorphisms
depending smoothly on $t\in\mathbb{R}$ (for example $\{\psi_{t}
\}_{t\in\mathbb{R}}$ could be the flow of a Hamiltonian vector field) and let
$S:\mathbb{R}\times Q\times Q\rightarrow\mathbb{R}$ be the corresponding
generating functions associated to this family. We will say that $\psi_{t}$
{\bfseries\itshape transforms a vector field $X\in\mathfrak{X}(T^{\ast}Q)$ to
equilibrium} if $T\psi_{t}\left(  X\right)  =0$ for any $t\in\mathbb{R}$. For
example, if $X=X_{h}$ is the Hamiltonian vector field associated to a
Hamiltonian function $h\in C^{\infty}\left(  T^{\ast}Q\right)  $ and $\psi
_{t}$ transforms $X_{h}$ to equilibrium, then the integral curve $\gamma$ of
$X_{h}$ with initial condition $z$ is
\[
\gamma_{t}=\hat{\psi}^{-1}\left(  \psi_{0}\left(  z\right)  ,t\right)
\]
where $\hat{\psi}^{-1}$ is the inverse of the diffeomorphism $\hat{\psi
}:T^{\ast}Q\times\mathbb{R}\rightarrow T^{\ast}Q\times\mathbb{R}$ given by
$\left(  z,t\right)  \mapsto\left(  \psi_{t}\left(  z\right)  ,t\right)  $.
The main goal of the classical Hamilton-Jacobi theory in this context is proving that $\psi$
transforms $X_{h}$ to equilibrium if, roughly speaking, its generating
function $S$ satisfies the (deterministic) Hamilton-Jacobi
equation. As we deal with time-dependent transformations $\psi_{t}$
of the phase space, the  time-dependent Hamiltonian formalism is
more convenient.

\bigskip

\noindent\textbf{Time-dependent Hamiltonian systems.} Recall that, for
time-dependent Hamiltonian systems, the phase space $T^{\ast}Q$ is replaced
with the \textit{extended} phase space $\mathbb{R}\times T^{\ast}Q$. Given a
time-dependent Hamiltonian function $h\in C^{\infty}\left(  \mathbb{R}\times
T^{\ast}Q\right)  $, one introduces $\Omega_{h}\in\Omega^{2}\left(
\mathbb{R}\times T^{\ast}Q\right)  $ as $\Omega_{h}=\mathbf{d}h\wedge
\mathbf{d}t+\omega$, where $\omega\in\Omega^{2}\left(  T^{\ast}Q\right)  $ is
the canonical symplectic form and $t$ denotes the global time coordinate in
$\mathbb{R}$. Observe that $\Omega_{h}$ is exact, $\Omega_{h}=-\mathbf{d}%
\theta_{h}$, where $\theta_{h}=\theta-h\mathbf{d}t$ and $\theta$ is the
canonical Liouville one form on the cotangent bundle. Then, the Hamiltonian vector field
$X_{h}\in\mathfrak{X}\left(  \mathbb{R}\times T^{\ast}Q\right)  $ is
characterized by the two equations
\[
i_{X_{h}}\Omega_{h}=0,\text{ \ \ }T\pi_{\mathbb{R}}\left(  X_{h}\right)
=\frac{\partial}{\partial t},
\]
where $\pi_{\mathbb{R}}:\mathbb{R}\times T^{\ast}Q\rightarrow\mathbb{R}$ is
the projection onto the first factor.

Sometimes it is more convenient to encode time-dependent Hamiltonian systems as autonomous Hamiltonian systems on the symplectic manifold $E:=T^{\ast}\left(
\mathbb{R}\times Q\right)  =T^{\ast}\mathbb{R}\times T^{\ast}Q$:  let $\left(  t,u\right)  $ be
global coordinates for $T^{\ast}\mathbb{R}$, that is $u$ is the conjugate
momentum associated to the time $t$, and denote by $\pi_{\mathbb{R}\times
T^{\ast}Q}:T^{\ast}\mathbb{R}\times T^{\ast}Q\rightarrow\mathbb{R}\times
T^{\ast}Q$ the projection $\left(  \left(  t,u\right)  ,z\right)
\mapsto\left(  t,z\right)  $, with $z\in T^{\ast}Q$. It is straightforward to check
that the Hamiltonian vector field $X_{h^{\star}}$ associated to the function
$h^{\star}:=u+\pi_{\mathbb{R}\times T^{\ast}Q}^{\ast}\left(  h\right) \in  C^{\infty}(E)$ is
such that $T\pi_{\mathbb{R}\times T^{\ast}Q}\left(  X_{h^{\star}}\right)
=X_{h}$. In other words, any time-dependent Hamiltonian system may be
visualized as an autonomous Hamiltonian system by  replacing
$\mathbb{R}\times T^{\ast}Q$ by $E$ and $h$ by $h^{\star}$; the integral
curves of the original system $X _h $ are simply obtained form the integral curves
of the autonomous system $X_{h^{\star}} $ by dropping the additional degree of freedom $u$, which is irrelevant as far as the dynamical description of the system is concerned. The following
proposition deals with a time-dependent family of symplectomorphisms $\left\{
\psi_{t}\right\}  _{t\in\mathbb{R}}$ of $T^{\ast}Q$ in the enlarged phase
space $E$ and will be useful in order to transform time-dependent Hamiltonian
systems.

\begin{proposition}
\label{prop HJ1}Let $\{\psi_{t}\}_{t\in\mathbb{R}}$ be a family of
symplectomorphisms of $T^{\ast}Q$ and $S\in C^{\infty}\left(  \mathbb{R}\times
Q\times Q\right)  $ its generating function. Define%
\[%
\begin{array}
[c]{rrl}%
\bar{\psi}~~:E & \longrightarrow & E\\
\left(  t,u,z\right)   & \longmapsto & \left(  t,u,\psi_{t}(z)\right)  ,
\end{array}
\]
where $t\in\mathbb{R}$, $u\in\mathbb{R}$,\ $z\in T^{\ast}Q$, and%
\begin{equation}%
\begin{array}
[c]{rrl}%
J_{t}:T^{\ast}Q & \longrightarrow & Q\times Q\\
z & \longmapsto & \left(  \pi(z),\pi\left(  \psi_{t}(z)\right)  \right)  .
\end{array}
\label{HJ 20}%
\end{equation}
Then,
\begin{enumerate}
\item[\textbf{(i)}] $\omega_{E}=\bar{\psi}^{\ast}\left(  \omega_{E}\right)
+\mathbf{d}\left(  \frac{\partial S}{\partial t}\circ J\circ\pi_{\mathbb{R}%
\times T^{\ast}Q}\right)  \wedge\mathbf{d}t$, where $\omega_{E}$ denotes the
canonical symplectic two form of $E=T^{\ast}\left(  \mathbb{R}\times Q\right)
$.

\item[\textbf{(ii)}] $\bar{\psi}^{\ast}\left(  \omega_{E}\right)  $ is
non-degenerate and, for any $\alpha\in\Omega\left(  \mathbb{R}\times T^{\ast
}Q\right)  $ and any $h\in C^{\infty}\left(  \mathbb{R}\times T^{\ast
}Q\right)  $,%
\[
\mathbf{d}h^{\star}\left(  \omega_{E}^{\#}\circ\bar{\psi}^{\ast}\circ
\pi_{\mathbb{R}\times T^{\ast}Q}\left(  \alpha\right)  \right)  =\mathbf{d}%
\left(  h\circ\hat{\psi}^{-1}+\frac{\partial S}{\partial t}\circ J_{t}%
\circ\hat{\psi}^{-1}\right)  ^{\star}\left(  \omega_{E}^{\#}\circ
\pi_{\mathbb{R}\times T^{\ast}Q}\left(  \alpha\right)  \right)  \circ\bar
{\psi}%
\]

\end{enumerate}
\end{proposition}

\noindent\textbf{Proof.} \ \ \ \ \textbf{(i)} Let $\left(  \left(  t,u\right)
,\left(  q^{i},p_{i};i=1,...,n\right)  \right)  $ be local coordinates on a
suitable open neighborhood $U\subseteq E$. It is immediate to see from
(\ref{HJ 24}) that for any $z\in T^{\ast}Q$%
\[
p_{i}\left(  z\right)  =\frac{\partial S}{\partial q_{1}^{i}}\left(
t,J_{t}\left(  z\right)  \right)  \text{ \ \ and \ \ }p_{i}\left(  \psi
_{t}(z)\right)  =-\frac{\partial S}{\partial q_{2}^{i}}\left(  t,J_{t}\left(
z\right)  \right)  ,
\]
$i=1,...,n$ (see, for instance, (7.9.1) in~\cite{mechanics and symmetry}), which implies that the
canonical one-form $\theta_{E}:=u\mathbf{d}t+\sum_{i=1}^{n}p_{i}%
\mathbf{d}q^{i}$ locally equals%
\[
\bar{\psi}^{\ast}\left(  \theta_{E}\right)  +\mathbf{d}S\circ J\circ
\pi_{\mathbb{R}\times T^{\ast}Q}-\frac{\partial S}{\partial t}\circ J\circ
\pi_{\mathbb{R}\times T^{\ast}Q}\mathbf{d}t
\]
(see, for instance, (7.9.5) in~\cite{mechanics and symmetry}). Applying $-\mathbf{d}$ to this expression, the result follows.

\textbf{(ii)} By (i), $\left(  \bar{\psi}^{-1}\right)  ^{\ast}\omega
_{E}=\omega_{E}+\mathbf{d}\left(  \frac{\partial S}{\partial t}\circ J\circ
\pi_{\mathbb{R}\times T^{\ast}Q}\circ\bar{\psi}^{-1}\right)  \wedge
\mathbf{d}t$. In order to simplify our notation let $F:=\frac{\partial S}{\partial t}\circ
J_{t}\circ\hat{\psi}^{-1}$. Then, using
$\{\mathbf{d}t,\mathbf{d}u,\mathbf{d}q^{i},\mathbf{d}p_{i}\}_{i=1,...,n}$ and
$\{\frac{\partial}{\partial t},\frac{\partial}{\partial u},\frac{\partial
}{\partial q^{i}},\frac{\partial}{\partial p_{i}}\}_{i=1,...,n}$ as bases of
$T^\ast _{\bar{\psi}(m)}U$ and $T_{\bar{\psi}(m)}U$ respectively, we have the
relations
\begin{equation}%
\begin{array}
[c]{ll}%
\left(  \left(  \bar{\psi}^{-1}\right)  ^{\ast}\omega_{E}\right)  ^{\#}\left(
\mathbf{d}t\right)  =-\frac{\partial}{\partial u}, & \omega_{E}^{\#}\left(
\mathbf{d}t\right)  =-\frac{\partial}{\partial u},\\
\left(  \left(  \bar{\psi}^{-1}\right)  ^{\ast}\omega_{E}\right)  ^{\#}\left(
\mathbf{d}u\right)  =\frac{\partial}{\partial t}+\sum_{i=1}^{n}\left(
\frac{\partial F}{\partial p_{i}}\frac{\partial}{\partial q^{i}}
-\frac{\partial F}{\partial q^{i}}\frac{\partial}{\partial p_{i}}\right),   &
\omega_{E}^{\#}\left(  \mathbf{d}u\right)  =\frac{\partial}{\partial t},\\
\left(  \left(  \bar{\psi}^{-1}\right)  ^{\ast}\omega_{E}\right)  ^{\#}\left(
\mathbf{d}q^{i}\right)  =-\frac{\partial F}{\partial p_{i}}\frac{\partial
}{\partial u}-\frac{\partial}{\partial p_{i}}, & \omega_{E}^{\#}\left(
\mathbf{d}q^{i}\right)  =-\frac{\partial}{\partial p_{i}},\\
\left(  \left(  \bar{\psi}^{-1}\right)  ^{\ast}\omega_{E}\right)  ^{\#}\left(
\mathbf{d}p_{i}\right)  =\frac{\partial F}{\partial q^{i}}\frac{\partial
}{\partial u}+\frac{\partial}{\partial q^{i}}, & \omega_{E}^{\#}\left(
\mathbf{d}p_{i}\right)  =\frac{\partial}{\partial q^{i}},
\end{array}
\label{HJ 25}%
\end{equation}
which easily shows the non-degeneracy of $\left(  \bar{\psi}^{-1}\right)  ^{\ast
}\omega_{E}$.

Let now $g\in C^{\infty}\left(  \mathbb{R}\times T^{\ast}Q\right)  $,
$\alpha\in\Omega\left(  \mathbb{R}\times T^{\ast}Q\right)  $, and $g^{\star
}=u+\pi_{\mathbb{R}\times T^{\ast}Q}^{\ast}(g)$. Using (\ref{HJ 25}), it is straightforward to check that
\begin{equation}
\mathbf{d}g^{\star}\left[  \left(  \left(  \bar{\psi}^{-1}\right)  ^{\ast
}\omega_{E}\right)  ^{\#}\left(  \pi_{\mathbb{R}\times T^{\ast}Q}^{\ast
}(\alpha)\right)  \right]  =\mathbf{d}\left(  g+F\right)  ^{\star}\left[
\omega_{E}^{\#}\left(  \pi_{\mathbb{R}\times T^{\ast}Q}^{\ast}(\alpha)\right)
\right]  .\label{HJ 27}%
\end{equation}
Additionally, for any $m\in U\subseteq E$,
the following diagram commutes:
\begin{equation}%
\begin{array}
[c]{rcl}%
T_{m}^{\ast}E & \overset{\omega_{E}^{\#}(m)}{\longrightarrow} & T_{m}E\\
_{T_{m}^{\ast}\bar{\psi}}\uparrow &  & \downarrow_{T_{m}\bar{\psi}}\\
T_{\bar{\psi}(m)}^{\ast}E & \overset{\left(  \left(  \bar{\psi}^{-1}\right)
^{\ast}\omega_{E}\right)  ^{\#}(\bar{\psi}(m))}{\longrightarrow} &
T_{\bar{\psi}(m)}E.
\end{array}
\label{HJ 26}%
\end{equation}
Therefore, by (\ref{HJ 26}), for any $\beta\in\Omega\left(  E\right)  $ and
any $h\in C^{\infty}\left(  \mathbb{R}\times T^{\ast}Q\right)  $,%
\begin{align*}
\mathbf{d}h^{\star}\left[  \omega_{E}^{\#}\circ\bar{\psi}^{\ast}%
(\beta)\right]  \left(  m\right)   &  =\mathbf{d}h^{\star}\left(  m\right)
\left[  \omega_{E}^{\#}\left(  m\right)  \left[  T_{m}^{\ast}\bar{\psi}\left(
\beta\left(  \bar{\psi}(m)\right)  \right)  \right]  \right]  \\
&  =\mathbf{d}h^{\star}\left(  m\right)  \left[  T_{\bar{\psi}(m)}\bar{\psi
}^{-1}\left[  \left(  \left(  \bar{\psi}^{-1}\right)  ^{\ast}\omega
_{E}\right)  ^{\#}\left(  \bar{\psi}(m)\right)  \left[  \beta\left(  \bar
{\psi}(m)\right)  \right]  \right]  \right]  \\
&  =\mathbf{d}\left(  \left(  \bar{\psi}^{-1}\right)  ^{\ast}h^{\star}\right)
\left(  \bar{\psi}(m)\right)  \left[  \left(  \left(  \bar{\psi}^{-1}\right)
^{\ast}\omega_{E}\right)  ^{\#}\left(  \bar{\psi}(m)\right)  \left[
\beta\left(  \bar{\psi}(m)\right)  \right]  \right]  \\
& = \mathbf{d}\left(  \left(  \bar{\psi}^{-1}\right)  ^{\ast}h^{\star}\right)
\left(  \bar{\psi}(m)\right)  \left[  \left(  \left(
\bar{\psi}^{-1}\right)  ^{\ast}\omega_{E}\right)  ^{\#}\left(  \bar{\psi
}(m)\right)  \left[  \beta\left(  \bar{\psi}(m)\right)  \right]  \right]  .
\end{align*}
In addition, if $\beta$ is of the form $\pi_{\mathbb{R}\times T^{\ast}Q}%
^{\ast}(\alpha)$ for some $\alpha\in\Omega\left(  \mathbb{R}\times T^{\ast
}Q\right)  $, by (\ref{HJ 27}) with $g=(\hat{\psi}^{-1})^{\ast}h$ we have%
\[
\mathbf{d}h^{\star}\left[  \omega_{E}^{\#}\circ\bar{\psi}^{\ast}\circ
\pi_{\mathbb{R}\times T^{\ast}Q}^{\ast}(\alpha)\right]  \left(  m\right)
=\mathbf{d}\left(  \left(  \left(  \hat{\psi}^{-1}\right)  ^{\ast}h+F\right)
^{\star}\right)  \left(  \bar{\psi}(m)\right)  \left[  \left(  \omega_{E}%
^{\#}\circ\pi_{\mathbb{R}\times T^{\ast}Q}^{\ast}(\alpha)\right)  \left(
\bar{\psi}(m)\right)  \right]  .
\]
Since $F=\frac{\partial S}{\partial t}\circ J_{t}\circ\hat{\psi}^{-1}$, the expression in (ii)
follows. \ \ \ \ $\blacksquare\medskip$

\begin{proposition}
\label{prop HJ2}Let $h\in C^{\infty}\left(  \mathbb{R}\times T^{\ast}Q\right)
$. With the same notation as in Proposition \ref{prop HJ1}, a curve
$\gamma:\left[  0,T\right]  \rightarrow\mathbb{R}\times T^{\ast}Q$ is a
solution of the Hamiltonian system defined by $h$ if and only if, for any
family of symplectomorphisms $\{\psi_{t}\}_{t\in\mathbb{R}}$ of $T^{\ast}Q$,
the curve $\psi\circ\gamma:\left[  0,T\right]  \rightarrow\mathbb{R}\times
T^{\ast}Q$ such that $(\hat{\psi}\circ\gamma)\left(  t\right)  :=(t,\psi
_{t}\left(  \gamma(t)\right)  )$ is a solution of a Hamiltonian system with
Hamiltonian function
\begin{equation}
h^{\prime}=h\circ\hat{\psi}^{-1}+\frac{\partial S}{\partial t}\circ J\circ
\hat{\psi}^{-1}\label{HJ 21}%
\end{equation}
where $S\in C^{\infty}\left(  \mathbb{R}\times Q\times Q\right)  $ is the
generating function of $\left\{  \psi_{t}\right\}  _{t\in\mathbb{R}}$.
\end{proposition}

\noindent\textbf{Proof.} \ \ \ Let $\gamma:\left[  0,T\right]  \rightarrow
\mathbb{R}\times T^{\ast}Q$ be a solution of the time-dependent Hamiltonian
system defined by $h$. Let $\bar{\gamma}:\left[  0,T\right]  \rightarrow
E=T^{\ast}\left(  \mathbb{R}\times Q\right)  $ be the curve such that
$\gamma=\pi_{\mathbb{R}\times T^{\ast}Q}^{\ast}(\bar{\gamma})$ and $\dot
{u}=\frac{\partial h}{\partial t}\left(  \gamma\right)  $, $u$ being the
conjugate momenta of the time coordinate $t$. Then $\gamma$ is a solution of
the time-dependent Hamiltonian system defined by $h\in C^{\infty}\left(
\mathbb{R}\times T^{\ast}Q\right)  $ if and only if $\bar{\gamma}$ is a
solution of the autonomous Hamilton system on the phase space $E$ with
Hamiltonian function $h^{\star}=u+\pi_{\mathbb{R}\times T^{\ast}Q}^{\ast}(h)$.
By (\ref{HJ 28}), this means that for any $\beta\in\Omega\left(  E\right)  $,
\begin{equation}
\int_{\left.  \bar{\gamma}\right\vert _{[0,t]}}\beta=-\int_{0}^{t}%
\mathbf{d}h^{\star}\left(  \omega_{E}^{\#}(\beta)\right)  \circ\gamma
(s)ds\label{HJ 29}%
\end{equation}
for any $t\in\lbrack0,T]$. However, since we are not interested in the
evolution of $u$, the conjugate momentum of the time, verifying that $\gamma$
is a solution of the time-dependent Hamilton equations is equivalent to taking
any curve $\bar{\gamma}$ such that $\gamma=\pi_{\mathbb{R}\times T^{\ast}
Q}^{\ast}(\bar{\gamma})$ and checking that (\ref{HJ 29}) holds for any
differential form of the type $\pi_{\mathbb{R}\times T^{\ast}Q}^{\ast}\left(
\alpha\right)  $, $\alpha\in\Omega\left(  \mathbb{R}\times T^{\ast}Q\right)  $.

Let now $\{\psi_{t}\}_{t\in\mathbb{R}}$ be a time-dependent family of
symplectomorphisms of $T^{\ast}Q$ and consider $\hat{\psi}:\mathbb{R}\times
T^{\ast}Q\rightarrow\mathbb{R}\times T^{\ast}Q$ such that $\hat{\psi}\left(
t,z\right)  =\left(  t,\psi_{t}\left(  z\right)  \right)  $, $\left(
t,z\right)  \in\mathbb{R}\times T^{\ast}Q$, and $\bar{\psi}:E\rightarrow E$
such that $\bar{\psi}\left(  t,u,z\right)  =\left(  t,u,\psi_{t}\left(
z\right)  \right)  $ as in Proposition \ref{prop HJ1}. Let $\bar{\psi}%
\circ\bar{\gamma}:[0,T]\rightarrow E$ be defined as $(\bar{\psi}\circ
\bar{\gamma})(s):=\bar{\psi}_{s}(\bar{\gamma}(s))$. Then%
\begin{align*}
\int_{\left.  \bar{\psi}\circ\bar{\gamma}\right\vert _{[0,t]}}\pi
_{\mathbb{R}\times T^{\ast}Q}^{\ast}\left(  \alpha\right)   &  =\int_{\left.
\bar{\gamma}\right\vert _{[0,t]}}\bar{\psi}^{\ast}\left(  \pi_{\mathbb{R}%
\times T^{\ast}Q}^{\ast}\left(  \alpha\right)  \right)  =-\int_{0}%
^{t}\mathbf{d}h^{\star}\left(  \omega_{E}^{\#}\circ\bar{\psi}^{\ast}%
(\pi_{\mathbb{R}\times T^{\ast}Q}^{\ast}\left(  \alpha\right)  )\right)
\circ\bar{\gamma}(s)ds\\
&  =-\int_{0}^{t}\mathbf{d}\left(  h\circ\hat{\psi}^{-1}+\frac{\partial
S}{\partial t}\circ J\circ\hat{\psi}^{-1}\right)  ^{\star}\left(  \omega
_{E}^{\#}(\pi_{\mathbb{R}\times T^{\ast}Q}^{\ast}\left(  \alpha\right)
)\right)  \circ\left(  \bar{\psi}\circ\bar{\gamma}\right)  (s)ds
\end{align*}
where Proposition \ref{prop HJ1} (ii) have been used in the last equality.
Hence, we conclude that $\pi_{\mathbb{R}\times T^{\ast}Q}(\bar{\psi}\circ
\bar{\gamma})=\hat{\psi}\circ\gamma$ is a solution of the time-dependent
Hamiltonian system given by $(\hat{\psi}^{-1})^{\ast}(h+\frac{\partial
S}{\partial t}\circ J)$. The converse is left to the reader.
\ \ \ \ $\blacksquare\medskip$

The content of Proposition \ref{prop HJ2} can be restated as follows. Given
$h\in C^{\infty}\left(  \mathbb{R}\times T^{\ast}Q\right)  $ and a family of
symplectomorphisms $\{\psi_{t}\}_{t\in\mathbb{R}}$, there exists a smooth
function $h^{\prime}\in C^{\infty}\left(  \mathbb{R}\times T^{\ast}Q\right)  $
such that $\Omega_{h}=\hat{\psi}^{\ast}\Omega_{h^{\prime}}$, where
$\Omega_{h^{\prime}}=\mathbf{d}h^{\prime}\wedge\mathbf{d}t+\omega$ and
$h^{\prime}$ is given by (\ref{HJ 21}) (see \cite[Section 7.9]{mechanics and
symmetry}). Furthermore $T\hat{\psi}^{-1}\left(  X_{h}\right)  $ is the
Hamiltonian vector field related to $h^{\prime}$ and the flow of
$X_{h^{\prime}}$ restricted to the phase space $T^{\ast}Q$ is $\hat{\varphi
}_{t}=\psi_{t}^{-1}\circ\varphi_{t}\circ\psi_{0}$ where, as usual, $\varphi$
denotes the flow of symplectomorphisms of the Hamiltonian vector field
$X_{h}\in\mathfrak{X}\left(  T^{\ast}Q\right)  $. However, as we will be
interested in transforming $X_{h}$ using $T\psi$ rather than $T\psi
^{-1}$ we will rewrite (\ref{HJ 21}) in the form
\begin{equation}
h^{\prime}\left(  t,\psi_{t}(z)\right)  :=h\left(  z\right)  +\frac{\partial
S}{\partial t}\left(  t,J_{t}\circ z\right)  .\label{HJ 22}%
\end{equation}

\begin{definition}
\label{def K}Let $h\in C^{\infty}\left(  T^{\ast}Q\right)  $ be a Hamiltonian
function and let $\left(  q^{i},p_{i};i=1,...,n\right)  $ be local Darboux
coordinates on $T^{\ast}Q$. Regarding $h$ as a function of these coordinates,
we will say that the generating function $S:\mathbb{R}\times Q\times
Q\rightarrow\mathbb{R}$ satisfies the (deterministic) Hamilton-Jacobi equation
if the function $K:\mathbb{R}\times Q\times Q\rightarrow\mathbb{R}$
\begin{equation}
K_{t}\left(  q_{1},q_{2}\right)  :=h\left(  q_{1},\frac{\partial S}{\partial
q_{1}}\left(  t,q_{1},q_{2}\right)  \right)  +\frac{\partial S}{\partial
t}\left(  t,q_{1},q_{2}\right)  \text{, \ }\left(  q_{1},q_{2}\right)  \in
Q\times Q \label{HJ 19}%
\end{equation}
does not depend on the first entry $q_{1}\in Q$.
\end{definition}

\noindent Observe that in the right hand side of (\ref{HJ 19}) we have carried
out the substitution $(p_{1})_{i}=\frac{\partial S}{\partial q_{1}^{i}%
}(t,q_{1},q_{2})$, $i=1,...,n$. We could also write (\ref{HJ 19}) more
intrinsically as%
\[
h\left(  \mathbf{d}_{Q_{1}}S\left(  t,q_{1},q_{2}\right)  \right)
+\frac{\partial S}{\partial t}\left(  t,q_{1},q_{2}\right)
\]
where, for a fixed value $\left(  t,q_{2}\right)  \in\mathbb{R}\times Q$, we
consider $\mathbf{d}_{Q_{1}}S\left(  t,q_{1},q_{2}\right)  $ as an element in
$T_{q_{1}}^{\ast}Q$.

Notice that the map $J_{t}$ introduced in (\ref{HJ 20}) is a local diffeomorphism for
any $t\in\mathbb{R}$ because we required the projection $\tau$ defined in
(\ref{HJ 23}) to be a local diffeomorphism when restricted to the graph of
$\psi_{t}$. We may therefore (locally) write any $z\in T^{\ast}Q$ as $z=J_{t}%
^{-1}\left(  q_{1},q_{2}\right)  $ for some suitable $\left(  q_{1}%
,q_{2}\right)  \in Q\times Q$. The important point is that $J_{t}^{-1}\left(
q_{1},q_{2}\right)  =\mathbf{d}_{Q_{1}}S\left(  t,q_{1},q_{2}\right)  $
(\cite[(7.9.1)]{mechanics and symmetry}) and, consequently, the transformed
Hamiltonian $h^{\prime}$ in~(\ref{HJ 22}) can be seen as a function on
$\mathbb{R}\times Q\times Q$. Explicitly, if $\bar{z}=\psi_{t}\left(
z\right)  \in T^{\ast}Q$,%
\begin{equation}
\label{write with z}
h^{\prime}\left(  t,\bar{z}\right)  =h\left(  \mathbf{d}_{Q_{1}}S\left(
t,q_{1},q_{2}\right)  \right)  +\frac{\partial S}{\partial t}\left(
t,q_{1},q_{2}\right)  ,
\end{equation}
so $h^{\prime}\left(  t,\bar{z}\right)  $ equals the function $K_{t}\left(
q_{1},q_{2}\right)  $ introduced in Definition \ref{def K}. Suppose now that
$S:\mathbb{R}\times Q\times Q\rightarrow\mathbb{R}$ is a solution to the
Hamilton-Jacobi equation. In other words, $K_{t}\left(  q_{1},q_{2}\right)
\equiv K_{t}\left(  q_{2}\right)  $. Since $q_{2}=\pi\left(  \psi_{t}\left(
z\right)  \right)  $ is the base point in the configuration space of the
transformed point $\psi_{t}\left(  z\right)  $, $z\in T^{\ast}Q$, we conclude
that $h^{\prime}$ does not depend on the fiber coordinates. Hence, removing
the subindices, the Hamilton equations associated to the new Hamiltonian
$h^{\prime}$ are 
\[
\dot{q}^{i}=0,\text{ \ \ }\dot{p}_{i}=-\frac{\partial K}{\partial q^{i}%
}\left(  t,q\right)  ,\text{ \ \ }i=1,...,n,
\]
which are easily integrable. In particular, if $K$ is independent of both
$q_{1}$ and $q_{2}$, then $\psi_{t}$ transforms $X_{h}$ to equilibrium.

\subsection{The stochastic case\label{seccion final}}

We are now going to see that the classical Hamilton-Jacobi that we just outlined has a stochastic counterpart. More specifically, one may use a
time-dependent family of symplectomorphisms and their generating function to
transform a stochastic Hamiltonian system into another one in much the same
fashion as in the deterministic case. The strategy consists of finding and characterizing a
suitable generating function so that the new Hamiltonian system is easier
to solve.

Let $T^{\ast}Q$ be the cotangent bundle of the configuration space manifold $Q$
and let $\{h_{0},h_{1},...,h_{r}\}\subset C^{\infty}\left(  T^{\ast}Q\right)
$ be a family of functions. Take a $\mathbb{R}^{r+1}$-valued semimartingale
$X:\mathbb{R}_{+}\times\Omega\rightarrow\mathbb{R}^{r+1}$ such that
\begin{equation}
\label{x in components}
X=\left(X^{0},X^{1},...,X^{r}\right), \quad\text{with} \quad X^{0}=t \quad\text{ a.s.,}
\end{equation}
and consider the
stochastic Hamiltonian system on $T^{\ast}Q$ with Hamiltonian function
$h:=(h_{0},h_{1},...,h_{r})$ and stochastic component $X$. If we want to
remove the assumption that there is a Hamiltonian vector field, i.e.
$X_{h_{0}}$, playing the role of a deterministic drift, we may simply choose
$h_{0}=0$. 

Using an approach similar to the one in Section \ref{seccion HJ determinista}, we will work in the extended phase space $E:=T^{\ast}\left(
\mathbb{R}\times Q\right)  $. Indeed, it is easy to check that the
solution semimartingales  of the stochastic Hamiltonian system can be
obtained out of the solutions of the stochastic Hamiltonian system on $E$ with
Hamiltonian function $\bar{h}=(h_{0}^{\star},\pi_{\mathbb{R}\times T^{\ast}%
Q}^{\ast}(h_{1}),...,\pi_{\mathbb{R}\times T^{\ast}Q}^{\ast}(h_{r}))$ and
stochastic component $X$; notice that the functions $h_{0},h_{1},...,h_{r}$ have
already been considered as functions on $\mathbb{R}\times T^{\ast}Q$ instead of only
$T^{\ast}Q$. The solutions of the original system can be recovered by
composing the solutions of the Hamiltonian system on $E$ with $\pi
_{\mathbb{R}\times T^{\ast}Q}$. When instead of working on the space $E$ one uses directly $\mathbb{R}\times T^{\ast}Q$ instead of $T^{\ast}Q$ then a $T^{\ast}Q$-valued semimartingale $\Gamma$ is a solution of the
corresponding stochastic Hamiltonian system when for any $\alpha\in
\Omega\left(  T^{\ast}Q\right)  $,
\[
\int\left\langle \alpha,\delta\Gamma_{s}\right\rangle =-\int\left\langle
\mathbf{d}h\left(  \tau_{T^{\ast}Q}^{\ast}\circ\omega^{\#}(\alpha)\right)
\left(  s,\Gamma_{s}\right)  ,\delta X_{s}\right\rangle ,
\]
where $\tau_{T^{\ast}Q}:\mathbb{R}\times T^{\ast}Q\rightarrow T^{\ast}Q$ is
the canonical projection onto the second factor.

\begin{proposition}
Let $\left\{  \psi_{t}\right\}  _{t\in\mathbb{R}}$ be a time-dependent family
of symplectomorphisms  of
$T^{\ast}Q$ with generating function $S\in C^{\infty}\left(  \mathbb{R}\times
Q\times Q\right)  $. Consider $\hat{\psi}:\mathbb{R}\times T^{\ast
}Q\rightarrow\mathbb{R}\times T^{\ast}Q$ and $\bar{\psi}:E\rightarrow E$ the
natural diffeomorphisms extending $\psi$ to $\mathbb{R}\times T^{\ast}Q$ and
$E$ respectively. Then the semimartingale $\Gamma:\mathbb{R}_{+}\times\Omega\rightarrow T^{\ast
}Q$ is a solution of the Hamiltonian system with Hamiltonian function
$h:T^{\ast}Q\rightarrow\mathbb{R}^{r+1}$, $h=\left(
h_{0},h_{1},...,h_{r}\right)  $, and stochastic component
$X:\mathbb{R}_{+}\times\Omega\rightarrow\mathbb{R}^{r+1}$ as in~(\ref{x in components}),  if and only if $\psi\left(  \Gamma\right)  $
is a solution of the Hamiltonian system with Hamiltonian function $h^{\prime
}:\mathbb{R}\times T^{\ast}Q\rightarrow\mathbb{R}^{r+1}$ with components given by%
\begin{align}
h_{0}^{\prime} &  =\tau_{T^{\ast}Q}^{\ast}(h_{0})\circ\hat{\psi}^{-1}
+\frac{\partial S}{\partial t}\circ J\circ\hat{\psi}^{-1},\nonumber\\
h_{1}^{\prime} &  =\tau_{T^{\ast}Q}^{\ast}(h_{1})\circ\hat{\psi}
^{-1},\nonumber\\
&  \vdots\nonumber\\
h_{r}^{\prime} &  =\tau_{T^{\ast}Q}^{\ast}(h_{r})\circ\hat{\psi}%
^{-1}.\label{HJ 32}%
\end{align}
and stochastic component $X$.
\end{proposition}

\noindent\textbf{Proof.} \ \ \ \ Suppose that $\Gamma:\mathbb{R}_{+}%
\times\Omega\rightarrow T^{\ast}Q$ is a solution of the Hamiltonian system
with Hamiltonian function $h$ and stochastic component $X$ and let
$\bar{\Gamma}:\mathbb{R}_{+}\times\Omega\rightarrow E$ be a semimartingale
such that $\pi_{\mathbb{R}\times T^{\ast}Q}(\bar{\Gamma}_{t})=(t,\Gamma
_{t})\in\mathbb{R}\times T^{\ast}Q$, $t\in\mathbb{R}_{+}$. We want to check
that $\bar{\psi}(\bar{\Gamma})$ is a solution of the stochastic Hamiltonian
system given by the Hamiltonian function (\ref{HJ 32}). Let $\alpha\in
\Omega\left(  \mathbb{R}\times T^{\ast}Q\right)  $. Since $\Gamma$ is a
solution, we may write%
\begin{gather*}
\int\left\langle \pi_{\mathbb{R}\times T^{\ast}Q}^{\ast}(\alpha),\delta
\bar{\psi}(\bar{\Gamma})\right\rangle =\int\left\langle \bar{\psi}^{\ast}%
\circ\pi_{\mathbb{R}\times T^{\ast}Q}^{\ast}(\alpha),\delta\bar{\Gamma
}\right\rangle =-\int\left\langle \mathbf{d}\bar{h}(\omega_{E}^{\#}\circ
\bar{\psi}^{\ast}\circ\pi_{\mathbb{R}\times T^{\ast}Q}^{\ast}(\alpha
))(\bar{\Gamma}),\delta X\right\rangle \\
=-\int\mathbf{d}h_{0}^{\star}(\omega_{E}^{\#}\circ\bar{\psi}^{\ast}\circ
\pi_{\mathbb{R}\times T^{\ast}Q}^{\ast}(\alpha))(\bar{\Gamma})dt-\sum
_{i=1}^{r}\int\mathbf{d}(\pi_{T^{\ast}Q}^{\ast}h_{i})(\omega_{E}^{\#}\circ
\bar{\psi}^{\ast}\circ\pi_{\mathbb{R}\times T^{\ast}Q}^{\ast}(\alpha
))(\bar{\Gamma})\delta X^{i},
\end{gather*}
where $\pi_{T^{\ast}Q}:E=T^{\ast}\mathbb{R}\times T^{\ast}Q\rightarrow
T^{\ast}Q$ is the projection onto the second factor. Now, by Proposition
\ref{prop HJ2} we have
\begin{multline}
\int\mathbf{d}h_{0}^{\star}(\omega_{E}^{\#}\circ\bar{\psi}^{\ast}\circ
\pi_{\mathbb{R}\times T^{\ast}Q}^{\ast}(\alpha))(\bar{\Gamma})dt=\\
\int\mathbf{d}\left(  \tau_{T^{\ast}Q}(h_{0})\circ\hat{\psi}^{-1}
+\frac{\partial S}{\partial t}\circ J_{t}\circ\hat{\psi}^{-1}\right)  ^{\star
}\left(  \omega_{E}^{\#}\circ\pi_{\mathbb{R}\times T^{\ast}Q}\left(
\alpha\right)  \right)  \left(  \bar{\psi}(\bar{\Gamma})\right)
dt.\label{HJ 33}%
\end{multline}
On the other hand, using (\ref{HJ 25}) and (\ref{HJ 26}) it is easy to
see that for any $g\in C^{\infty}\left(  T^{\ast}Q\right)  $
\[
\mathbf{d}(\pi_{T^{\ast}Q}^{\ast}g)(\omega_{E}^{\#}\circ\bar{\psi}^{\ast}%
\circ\pi_{\mathbb{R}\times T^{\ast}Q}^{\ast}(\alpha))(m)=\mathbf{d}%
(\pi_{T^{\ast}Q}^{\ast}(g)\circ\bar{\psi}^{-1})(\omega_{E}^{\#}\circ
\pi_{\mathbb{R}\times T^{\ast}Q}^{\ast}(\alpha))(\bar{\psi}(m)).
\]
Consequently,
\begin{equation}
\int\mathbf{d}(\pi_{T^{\ast}Q}^{\ast}h_{i})(\omega_{E}^{\#}\circ\bar{\psi
}^{\ast}\circ\pi_{\mathbb{R}\times T^{\ast}Q}^{\ast}(\alpha))(\bar{\Gamma
})\delta X^{i}=\int\mathbf{d}(\pi_{T^{\ast}Q}^{\ast}(h_{i})\circ\bar{\psi
}^{-1})(\omega_{E}^{\#}\circ\pi_{\mathbb{R}\times T^{\ast}Q}^{\ast}%
(\alpha))(\bar{\psi}(\bar{\Gamma}))\delta X^{i},\label{HJ 34}
\end{equation}
for any $i=1,...,r$. Combining (\ref{HJ 33}) and (\ref{HJ 34}) we
obtain that%
\begin{align*}
\int\left\langle \pi_{\mathbb{R}\times T^{\ast}Q}^{\ast}(\alpha),\delta
\bar{\psi}(\bar{\Gamma})\right\rangle  &  =-\int\mathbf{d}\left(
\tau_{T^{\ast}Q}(h_{0})\circ\hat{\psi}^{-1}+\frac{\partial S}{\partial t}\circ
J_{t}\circ\hat{\psi}^{-1}\right)  ^{\star}\left(  \omega_{E}^{\#}\circ
\pi_{\mathbb{R}\times T^{\ast}Q}\left(  \alpha\right)  \right)  \left(
\bar{\psi}(\bar{\Gamma})\right)  dt\\
&  -\int\mathbf{d}(\pi_{T^{\ast}Q}^{\ast}(h_{i})\circ\bar{\psi}^{-1}%
)(\omega_{E}^{\#}\circ\pi_{\mathbb{R}\times T^{\ast}Q}^{\ast}(\alpha
))(\bar{\psi}(\bar{\Gamma}))\delta X^{i},
\end{align*}
which means that $\psi_{t}\left(  \Gamma_{t}\right)  $ is a solution of the
time-dependent stochastic Hamiltonian system with stochastic component $X$ and
Hamiltonian function (\ref{HJ 32}). The converse is left to the reader.
\quad $\blacksquare$

\medskip

The system (\ref{HJ 32}) may be written as%
\begin{align}
h_{0}^{\prime}(t,\psi_{t}(z)) &  =\tau_{T^{\ast}Q}^{\ast}(h_{0})\left(
\mathbf{d}_{Q_{1}}S\left(  t,q_{1},q_{2}\right)  \right)  +\frac{\partial
S}{\partial t}\left(  t,q_{1},q_{2}\right)  \nonumber\\
h_{1}^{\prime}(t,\psi_{t}(z)) &  =\tau_{T^{\ast}Q}^{\ast}(h_{1})\left(
\mathbf{d}_{Q_{1}}S\left(  t,q_{1},q_{2}\right)  \right)  \nonumber\\
&  \vdots\nonumber\\
h_{r}^{\prime}(t,\psi_{t}(z)) &  =\tau_{T^{\ast}Q}^{\ast}(h_{r})\left(
\mathbf{d}_{Q_{1}}S\left(  t,q_{1},q_{2}\right)  \right)  \label{HJ 35}%
\end{align}
where, as in (\ref{write with z}) we have written
$z\in T^{\ast}Q$ as $z=J_{t}^{-1}\left(  q_{1},q_{2}\right)  $ for some
suitable $\left(  q_{1},q_{2}\right)  \in Q\times Q$. In addition, if the
generating function $S$ is such that the right hand side of (\ref{HJ 35}) is
independent of the variable $q_{1}$, that is,
\begin{align}
\tau_{T^{\ast}Q}^{\ast}(h_{0})\left(  \mathbf{d}_{Q_{1}}S\left(  t,q_{1}%
,q_{2}\right)  \right)  +\frac{\partial S}{\partial t}\left(  t,q_{1}%
,q_{2}\right)   &  =:K_{0}\left(  t,q_{2}\right),  \nonumber\\
\tau_{T^{\ast}Q}^{\ast}(h_{1})\left(  \mathbf{d}_{Q_{1}}S\left(  t,q_{1}%
,q_{2}\right)  \right)   &  =:K_{1}\left(  t,q_{2}\right),  \nonumber\\
&  \vdots\nonumber\\
\tau_{T^{\ast}Q}^{\ast}(h_{r})\left(  \mathbf{d}_{Q_{1}}S\left(  t,q_{1}%
,q_{2}\right)  \right)   &  =:K_{r}\left(  t,q_{2}\right)  ,\label{HJ 37}%
\end{align}
then the stochastic Hamilton equations of the transformed system may be expressed
in local coordinates as%
\begin{align*}
\delta q^i &  =0\\
\delta p_{i} &  =-\frac{\partial K_{0}}{\partial q}\left(  t,q\right)
dt-\sum_{i=1}^{r}\frac{\partial K_{i}}{\partial q}\left(  t,q\right)  \delta
X^{i}.
\end{align*}

\medskip

\noindent The next result is basically due to Bismut (see {\cite[Th\'{e}or\`{e}me 7.6, page 349]{Bismut}}).

\begin{proposition}
In the conditions of the previous proposition, if
(\ref{HJ 37}) holds then
\begin{align*}
\left\{  h_{i},h_{j}\right\}  (z)  &  =0\\
\{h_{0},h_{i}\}(z)+\frac{\partial K_{i}}{\partial t}\left(  t,\pi\left(
\psi_{t}(z)\right)  \right)   &  =0
\end{align*}
locally for any $1\leq i,j\leq r.$
\end{proposition}

\noindent\textbf{Proof.} \ \ \ \ Suppose that there exists a generating
function $S\in C^{\infty}\left(  \mathbb{R}\times Q\times Q\right)  $ such
that the equalities (\ref{HJ 37}) are satisfied. We take a fixed point
$q_{2}\in Q$ and write $K_{i}^{q_{2}}\left(  t\right)  $ instead of
$K_{i}\left(  t,q_{2}\right)  $, $i=0,...,r$, and $S^{q_{2}}\left(
t,q\right)  $ instead of $S\left(  t,q,q_{2}\right)  $. Consider the following
family of functions of the extended phase space $E=T^{\ast}\left(
\mathbb{R}\times Q\right)  $:%
\begin{align*}
g_{0} &  =u+\pi_{T^{\ast}Q}^{\ast}(h_{0})-K_{0}^{q_{2}}\left(  t\right)  \\
g_{1} &  =\pi_{T^{\ast}Q}^{\ast}(h_{1})-K_{1}^{q_{2}}\left(  t\right)  \\
&  \vdots\\
g_{r} &  =\pi_{T^{\ast}Q}^{\ast}(h_{r})-K_{r}^{q_{2}}\left(  t\right)  ,
\end{align*}
where $u$ denotes the conjugate momentum of the time coordinate $t$ in $E$. The functions $g_{0},...,g_{r}\subset C^{\infty}(E)$ vanish
on the Lagrangian submanifold $L_{S}\subset E$ locally defined by%
\[
L_{S}=\left\{  \left(  t,u,q,p\right)  \in E~|~p_{i}=\frac{\partial S^{q_{2}}%
}{\partial q^{i}}(t,q),~u=\frac{\partial S^{q_{2}}}{\partial t}(t,q)\right\}
.
\]
Given that if a family of functions is locally constant
on a Lagrangian submanifold, then their Poisson brackets must vanish on
it, we have that  $\{g_{i},g_{j}\}=0$ for any $0\leq i,j\leq r$. Equivalently,
\begin{align}
0 &  =\left.  \{\pi_{T^{\ast}Q}^{\ast}h_{i},\pi_{T^{\ast}Q}^{\ast}%
h_{j}\}\right\vert _{L_{S}}=\left.  \pi_{T^{\ast}Q}^{\ast}\left(  \left\{
h_{i},h_{j}\right\}  \right)  \right\vert _{L_{S}},\nonumber\\
0 &  =\left.  \pi_{T^{\ast}Q}^{\ast}\left(  \left\{  h_{0},h_{i}\right\}
\right)  \right\vert _{L_{S}}+\left.  \frac{\partial K_{i}^{q_{2}}}{\partial
t}\right\vert _{L_{S}},\label{HJ 38}%
\end{align}
for any $i,j=1,...,r$. In particular, since the inverse $J_{t}^{-1}:Q\times
Q\rightarrow T^{\ast}Q$ of the local diffeomorphism introduced in
(\ref{HJ 20}) is such that $z=J_{t}^{-1}\left(  q_{1},q_{2}\right)  =\left(
q_{1},\mathbf{d}S^{q_{2}}\left(  t,q_{1}\right)  \right)  $, we have the
freedom to chose $q_{2}$ so that $z=J_{t}^{-1}\left(  q_{1},q_{2}\right)  $
is a point in the fiber of $q_{1}\in Q$. With this choice (\ref{HJ 38})
implies that%
\begin{align*}
\left\{  h_{i},h_{j}\right\}  \left(  z\right)   &  =0\\
\{h_{0},h_{i}\}\left(  z\right)  +\frac{\partial K_{i}^{q_{2}}}{\partial
t}\left(  t\right)   &  =\{h_{0},h_{j}\}\left(  z\right)  +\frac{\partial
K_{i}}{\partial t}\left(  t,\pi\left(  \psi_{t}(z)\right)  \right)  =0
\end{align*}
for any $z\in T^{\ast}Q$. \ \ \ \ $\blacksquare$

\end{document}